\documentclass{article}
\usepackage{amsmath}
\usepackage{amssymb}
\usepackage{amscd}
\usepackage{multicol}
\usepackage{color}
\begin{document}
\def\ni{\noindent}
\def\t{\theta}
\def\O{\Omega}
\def\S{\Sigma}
\def\e{\epsilon}
\def\lra{\longrightarrow}
\def\R{{\mathbb R}}
\def\N{{\mathbb N}}
\def\Z{{\mathbb Z}}
\def\C{{\mathbb C}}
\def\RR{{{\mathbb R}}^2}
\def\MS{M^{2{\rm x}2}_s}
\def\N{{\bf N}}
\def\l{\lambda}
\def\LL{${\cal L}$}
\def\E{{\cal E}}
\def\a{{\alpha}}
\def\A{A_{\a}}
\def\ta{\t^{\a}}
\def\rot{{\rm rot}}

\newcommand{\vs}[1]{\vskip #1pt}
\newcommand\qed{\hfill \rule{6pt}{6pt}}

\newtheorem{Theorem}{Theorem}[section]
\newtheorem{Definition}[Theorem]{Definition}
\newtheorem{corollary}[Theorem]{Corollary}
\newtheorem{proposition}[Theorem]{Proposition}
\newtheorem{examples}[Theorem]{Esempi}
\newtheorem{example}[Theorem]{Example}
\newtheorem{lemma}[Theorem]{Lemma}
\newtheorem{remark}[Theorem]{Remark}

\catcode`\@=11


   \renewcommand{\theequation}{\thesection.\arabic{equation}}
   \renewcommand{\section}%
   {\setcounter{equation}{0}\@startsection {section}{1}{\z@}{-3.5ex
plus -1ex
    minus -.2ex}{2.3ex plus .2ex}{\Large\bf}}

\title{Linear and nonlinear eigenvalue problems for Dirac systems in unbounded domains\thanks{AMS Subject Classification: 34C23, 34B09, 34L40.}
\thanks{Under the auspices of GNAMPA-I.N.d.A.M., Italy.
The work of the first two authors has been performed in the frame of the M.I.U.R. Project
'Topological and Variational Methods in the Study of Nonlinear Phenomena'; the work of the third author has been performed in the frame of the M.I.U.R. Project 'Nonlinear Control: Geometrical Methods and Applications' and of the GNAMPA-I.N.d.A.M. project ``Equazione di evoluzione degeneri e singolari: controllo e applicazioni''.
}
}
\author{Anna Capietto$\ddag$, Walter Dambrosio$\ddag$ and Duccio Papini** \\
\\
$\ddag$ {\small Dipartimento di Matematica - Universit\`a di Torino}\\
{\small Via Carlo Alberto 10, 10123 Torino - Italia}\\
{\small E-mail addresses: anna.capietto@unito.it,  
walter.dambrosio@unito.it}
\\
**{\small Dipartimento di Ingegneria dell'Informazione e Scienze Matematiche - Universit\`a di
Siena}\\
{\small Via Roma 56, 53100 Siena - Italia}\\
{\small E-mail address: duccio.papini@unisi.it}}
\date{}
\maketitle


\vs{12}

\centerline{\bf Abstract}
\vs{12}
\ni
{\small
We first study the linear eigenvalue problem for a planar Dirac system in the open half-line and describe the nodal properties of its solution by means of the rotation number. We then give a global bifurcation result for a planar nonlinear Dirac system in the open half-line. As an application, we provide a global continuum of solutions of the nonlinear Dirac equation which have a special form.
}
\vs{12}
\ni {\bf Keywords.} Dirac system, eigenvalue problem, rotation number, global bifurcation.
\vs{12}
\ni

\vs{12}
\ni
\section{Introduction}
\vs{12}
\ni
In this paper we give a global bifurcation result (Theorem \ref{main-alternative}) for a nonlinear Dirac system in ${\R}^2$ of
the form
\begin{equation} \label{Ieqcompleta}
Jz'+P(x)z=\l z+S(x,z)z,\quad x>0,\quad \l\in \R,\quad z=
(u,v)\in {\R}^2,
\end{equation}
where
\[
J=\left(\begin{array}{cc}
              0&1\\
             &\\
             -1&0
            \end{array}
     \right)
\]
and $P(x), S(x,z)$ are
continuous symmetric matrices, for every $x>0$ and $z\in {\R}^2$. We will be interested  
in solutions $z$ of \eqref{Ieqcompleta} belonging to the space
\[
D_0 = \{z\in L^2(0,+\infty):\ z\in AC(0,+\infty),\ J z'+P(\cdot)z \in L^2(0,+\infty)\}.
\]
\ni In particular, the solutions are convergent to zero at zero and at
infinity. This choice is strictly related to the spectral properties  
of the linear operator $\tau z=Jz'+P(x)z$ and to the possibility of  
considering self-adjoint extensions of $\tau$ (see Section \ref{sez-autovalori}).

\ni When $P$ has the form
\begin{equation} \label{Dirac-magn_in}
P(x)=P_{V,k,\mu_a} (x)=\left(\begin{array}{cc}
 	      -1+V(x)&{\displaystyle{-\dfrac{k}{x}-\mu_a V'(x)}}\\
&\\
			         {\displaystyle{-\dfrac{k}{x}}-\mu_a V'(x)}&1+V(x)
				\end{array}
			\right),\quad x>0,
\end{equation}
the differential operator $z \mapsto Jz'+P(\cdot)z$ coincides with the radially symmetric Dirac operator with or without anomalous magnetic moment (cf. \cite{KaSc-04,ScTr-02, Th-book-92,We-book-87} and Section \ref{Dirac-nonlineare}). In this context $V\in C^1(0,+\infty)$ represents an electrostatic potential, $\mu_a\in {\R}$ an anomalous magnetic moment and $k\in {\Z}\setminus \{0\}$ (see \cite{Th-book-92}). For a comprehensive treatment of linear and nonlinear Dirac systems, we refer to the paper by M. Esteban \cite{Es-02}. As it is explained in detail in Section \ref{Dirac-nonlineare}, nonlinear systems of the form \eqref{Ieqcompleta} arise, for some $S$, when one is interested in solutions of a nonlinear Dirac PDE which have a special form (cf. \eqref{conto-nonlinearita-4}).

\ni The study of global bifurcation problems for second order equations in unbounded intervals was initiated in the 70s by C. Stuart \cite{St-75} and N. Dancer \cite{Da-75,Da-77}. More recent results have been given by P. Rabier-C. Stuart \cite{RaSt-01}, S. Secchi-C. Stuart \cite{SeSt-03}, the first and second author \cite{CaDa-10} and the authors \cite{CaDaPa-12}.

\ni In \cite{CaDa-10} it is considered the particular case when the r.h.s. of \eqref{Ieqcompleta} (and the function $S$) is regular at zero. We are now able to avoid this restriction and, as a consequence, to treat the physically relevant Dirac operator.

\ni Having in mind a bifurcation result, a comprehensive knowledge of the linear eigenvalue problem 
\begin{equation} \label{lineare1_in}
Jz'+P(x)z=\l z,\quad x>0,\quad \l\in \R,\quad z=(u,v)\in {\R}^2
\end{equation}
is necessary. More precisely, we have to study the existence of eigenvalues and their "nodal properties". To this end, in Subsection \ref{stimeasintotiche}, assuming $({\cal P}_1),({\cal P}_2),({\cal P}_3)$ for the matrix $P$,  we first describe (Lemmas \ref{limitiuv} and \ref{limitiuvzero}) the behaviour of the solutions of the linear system \eqref{lineare1_in} when $x\to +\infty$ and $x\to 0^+$. As in \cite{CaDa-10}, we apply the Levinson theorem \cite{Ea-book-89} on the asymptotic properties of solutions of linear equations and, by means of a suitable change of variables, we manage to treat the singularity at zero as well. Using the results of Subsection \ref{stimeasintotiche}, we develop in Subsection \ref{oscillazione} an oscillatory theory for nontrivial solutions of \eqref{lineare1} based on the study of the asymptotic behaviour of the angular coordinate $\theta$ in the phase-plane (cf. the book by J. Weidmann \cite{We-book-87}). It is interesting to observe that, contrary to the case of second order equations, 
in case of planar Dirac-type systems the angular coordinate is not, in general, an increasing function of $x$. However, we are able to guarantee (Propositions \ref{limitethetainf} and \ref{limitethetazero}) that the limits
\begin{equation}
\theta (+\infty,\l)=\lim_{x\to +\infty} \theta(x,\l),  \quad \theta (0)=\lim_{x\to 0^+} \theta(x,\l)
\end{equation}
exist and are finite. We can thus give the definition of  
\begin{equation} \label{defrotazioni_in}
\rot \ (z)=\dfrac{\theta(+\infty,\l)-\theta(0)}{\pi},
\end{equation}
the rotation number of a solution $z$ to $\eqref{lineare1_in}$. Roughly speaking, the unboundedness of the interval and the singularity at zero do not prevent solutions to perform only a finite number of rotations around the origin (as in the regular case). A nontrivial phase-plane analysis leads then to some useful continuity properties of the angular function near zero and infinity (Propositions \ref{confini-theta-inf} and \ref{confini-theta-zero}).

\ni In Section \ref{sez-autovalori} we study the spectral theory for the linear operator formally defined by
\begin{equation} \label{tau_in}
\tau z=Jz'+P(x)z,\ x>0.
\end{equation}
\ni More precisely, standard arguments from \cite{We-book-87} ensure that $\tau$ is in the limit point case at infinity and at zero and that there exists a unique self-adjoint realization $A_0$ (cf. \eqref{def-A0}) of $\tau$ having (when $P$ has the form \eqref{Dirac-magn_in}) essential spectrum $\sigma_{\rm{ess}} (A_0)=(-\infty,-1]\cup [1,+\infty)$. Then, the (nontrivial) question of characterizing eigenvalues of $A_0$ is tackled by the results of Subsection \ref{oscillazione}. Finally, we give results on the existence and accumulation of eigenvalues of $A_0$ at the boundary of the interval $(-1,1)$ which are based on the oscillatory behaviour of the solutions for a value of $\l$ corresponding to one of the extrema of the essential spectrum; similar results can be found in the case of second-order differential operators in the book by N. Dunford-J. Schwartz \cite{DuSc-book-63} and in case of Dirac operators (without any knowledge of the nodal properties of the corresponding eigenvalues) in the paper by H. Schmid-C. Tretter \cite{ScTr-02}.

\ni Taking advantage of all the results described above, in Subsection 4.1 we give a global bifurcation result (Theorem \ref{main-alternative})  for system \eqref{Ieqcompleta}. Due to the fact that we are dealing with an unbounded interval, we face a lack of  
compactness; this difficulty is overcome by applying an abstract bifurcation result due to C. Stuart \cite{St-75}. A more precise description of the continuum emanating from eigenvalues of odd multiplicity of the linear operator $\tau$ is then performed (as we did in \cite{CaDa-10}) in Theorem \ref{main-alternative-2}; to this aim, we develop a continuity-connectivity argument based on a linearization approach and on the properties of the rotation number of a solution to \eqref{Ieqcompleta} (cf. \eqref{def-rot-2}, \eqref{def-i-1} and Proposition \ref{continuita-i}).

\ni Finally, in Subsection 4.2 we consider the partial differential equation 
\begin{equation} \label{Dirac-1_in}
i\sum_{j=1}^3 \alpha_j \dfrac{\partial \psi}{\partial x_j}-\beta \psi - V(||x||)\psi + i a \sum_{j=1}^3 \alpha_j \dfrac{\partial V(||x||)}{\partial x_j} \psi  =\l \psi+\gamma (||x||)F(\langle\beta\psi,\psi\rangle)\beta\psi,\quad x\in {\R}^3, \ a\in \R,
\end{equation}
where $\psi:{\R}^3\to {\C}^4$, $V\in C((0,+\infty),\R)$ and $\gamma \in C((0,+\infty),\R)$ satisfy suitable assumptions, $\langle\cdot,\cdot\rangle$ denotes the scalar product in ${\C}^4$ and $\alpha_j$ ($j=1, 2, 3$) and $\beta$ are the $4\times 4$ Dirac matrices (see Subsection 4.2). Set 
\begin{equation} \label{def-H0_in}
H_0\psi=i\sum_{j=1}^3 \alpha_j \dfrac{\partial \psi}{\partial x_j}-\beta \psi,\quad \forall \ \psi \in H^1_0({\R}^3)\subset L^2({\R}^3).
\end{equation}
\ni It is well-known (cf. the book by B. Thaller \cite{Th-book-92}), that there exist suitable subspaces of $L^2(S^2)$ s.t. the restriction of the linear operator $H_0 -V+ia\ \alpha\cdot \nabla V $ to each of these subspaces can be represented by an ordinary differential operator of the form $\tau$. A remark on the physical meaning of the partial wave subspaces can be found in Remark \ref{significato} in Section \ref{Dirac-nonlineare}. It is interesting to observe (on the lines of a paper by F. Cacciafesta \cite{Ca-11}) that there are {\it nonlinear} terms $F(\langle\beta\psi,\psi\rangle)\beta\psi$ in \eqref{Dirac-1_in} which leave the above described subspaces invariant. These appear, among others, in the so-called Soler model and are the most interesting from a physical point of view (cf. \cite{Ra-book-83},\cite{So-70}). On the same lines, we refer also to the contributions by M. Balabane-T. Cazenave-L. Vazquez \cite{BaCaVa-90}, Y. Ding-B. Ruf \cite{DiRu-08}, J. Ding, J. Xu, F. Zhang \cite{DXZ-09}, Y. Dong-J. Xie \cite{DoXi-11} and references therein.

\noindent Our contribution (Theorem \ref{bif-Dirac}) provides the existence of a global continuum of solutions of the nonlinear PDE \eqref{Dirac-1_in} which have a special form (i.e. which belong to one of the above mentioned subspaces). To the authors' knowledge, Theorem \ref{bif-Dirac} is the first global bifurcation result for a nonlinear Dirac-type equation of the form \eqref{Dirac-1_in}. In the particular case $V \equiv 0$, M. Balabane-T. Cazenave-A. Douady-F. Merle \cite{BCDM-88} gave a multiplicity result for solutions (having prescribed nodal properties) to a system of ODEs of the form \eqref{eqcompleta-Dirac}. For multiplicity results via critical point theory for the nonlinear Dirac PDE, we refer to Theorem 3.3 in \cite{Es-02} (in case $V \equiv 0$) and to the paper by Y. Ding-B. Ruf \cite{DiRu-08} (for a potential that includes the Coulomb case). On the other hand, in the particular case of linear Dirac-type systems of ODEs, H. Schmid-C. Tretter \cite{ScTr-02} have given results for the eigenvalue problem for some special choice of the potential $V$. 

\vs{12}
\ni
In what follows, we will denote by $M^2_S$ the set of symmetric $2
\times 2$ matrices.

\section{Linear Dirac systems} \label{sistemalineare}
\vs{12}
\ni
In this Section we consider a linear system of the form
\begin{equation} \label{lineare1}
Jz'+P(x)z=\l z,\quad x>0,\quad \l\in \R,\quad z=(u,v)\in {\R}^2;
\end{equation}
by a solution of \eqref{lineare1} we mean a function $z\in AC_{{\rm loc}} (0,+\infty)$ satisfying \eqref{lineare1} almost everywhere in $(0,+\infty)$. In the next Sections we will be interested in solutions $z\in L^2(0,+\infty)$ or $z\in H^1(0,+\infty)$; hence, in describing the solutions of \eqref{lineare1} we will point out, when possible, if they belong to $L^2(0,+\infty)$ or to $H^1(0,+\infty)$.
\vs{6}
\ni
We assume that $P\in C((0,+\infty),M^2_S)$ and we denote by $p_{ij}$ its coefficients, as usual.
For each pair of real numbers $ \mu_{-} < \mu^{+} $, let us consider the class ${\cal P}_{\mu}$ of continuous maps $P:(0,+\infty)\lra M^{2,2}_S$ satisfying the following conditions:
\begin{description}
\item{$({\cal P}_1)$} There exists $q_\infty\geq 1$ such that
\begin{equation} \label{limitiPinfinito}
\lim_{x\to +\infty} P(x)=\left(\begin{array}{cc}
						\mu^-&0\\
						&\\
						0&\mu^+
						\end{array}
			\right)=:P_\infty
\end{equation}
and
\begin{equation} \label{iporestoinfinito}
\int_1^{+\infty} ||R_\infty (x)||^{q_\infty}\,dx<+\infty,
\end{equation}
where $R_\infty(x)=P(x)-P_\infty$, for every $x\geq 1$.
\item{$({\cal P}_2)$} There exist $\beta\geq 1$, $P^*\in M^2_S$ and $q_0\geq 1$ such that
\begin{equation} \label{limitiPzero}
\lim_{x\to 0^+} x^\beta  P(x)=P^*
\end{equation}
and
\begin{equation} \label{iporestozero}
\int_0^{1} \dfrac{1}{x^\beta} ||R_0(x)||^{q_0}\,dx<+\infty,
\end{equation}
where $R_0(x)=x^{\beta} P(x)-P^*$, for every $x\in (0,1)$.
\item{$({\cal P}_3)$} The matrix $P^*$ satisfies
\begin{eqnarray}
{\rm det }\  P^*<-1/4 \quad {\mbox{if $\beta=1$}} \label{determinante-beta-1} \\
\nonumber \\
{\rm det }\  P^*<0\quad {\mbox{if $\beta>1.$}} \label{determinante-beta}
\end{eqnarray}
\end{description}
\vs{12}
\ni
\noindent In what follows, we write $\Lambda=(\mu^-,\mu^+).$
\vs{12}
\ni
\begin{remark} \label{oss-ipotesi} 1. We observe that assumption $({\cal P}_2)$ implies that \eqref{lineare1} has a singularity for $x\to 0^+$; indeed, from \eqref{limitiPzero} and the fact that $P^*$ is not the zero-matrix (since its determinant is negative in any case), we deduce that
\[
p_{ij}(x)\sim \dfrac{p_{ij}^*}{x^\beta},\quad x\to 0^+ \qquad(i,j=1, 2)
\]
and, in particular, that $p_{ij}\notin L^1(0,1)$.
\vs{6}
\ni
(2) Let us also observe that, for a particular choice of $P$, the differential operator given in \eqref{lineare1} coincides with the radially symmetric Dirac operator with or without anomalous magnetic moment (cf. \cite{KaSc-04, ScTr-02, Th-book-92, We-book-87} and Section \ref{Dirac-nonlineare}); indeed, this is the situation when $P$ has the form
\begin{equation} \label{Dirac-magn}
P(x)=P_{V,k,\mu_a} (x)=\left(\begin{array}{cc}
 	      -1+V(x)&{\displaystyle{-\dfrac{k}{x}-\mu_a V'(x)}}\\
&\\
			         {\displaystyle{-\dfrac{k}{x}}-\mu_a V'(x)}&1+V(x)
				\end{array}
			\right),\quad x>0,
\end{equation}
where $V\in C^1(0,+\infty)$ is an electrostatic potential, $\mu_a\in {\R}$ is an anomalous magnetic moment and $k\in {\Z}\setminus \{0\}$ (see \cite{Th-book-92}).

\ni
Let us assume that $V$ satisfies the following conditions:
\begin{equation} \label{ipotesi-potenziale-infinito}
\begin{array}{l}
\displaystyle{V(x)=\dfrac{\gamma_\infty}{x^{\alpha_\infty}} +R_{V,\infty} (x),\ \alpha_\infty >0,}\\
\\
\displaystyle{x^{\alpha_\infty} R_{V,\infty}=o(1),\quad x^{\alpha_\infty+1} R'_{V,\infty}=o(1),\ x\to +\infty}
\end{array}
\end{equation}
and
\begin{equation} \label{ipotesi-potenziale-zero}
\begin{array}{l}
\displaystyle{V(x)=\dfrac{\gamma_0}{x^{\alpha_0}} +R_{V,0} (x),\ \alpha_0 >0,}
\end{array}
\end{equation}
where
\begin{equation} \label{ipotesi-potenziale-zero-resto1}
 {\mbox{if $\mu_a=0$}}: \quad \left\{\begin{array}{l}
\displaystyle{\alpha_0=1,\quad x R_{V,0}=o(1),\ x\to 0^+,\quad }\\
\\
\displaystyle{\int_0^1 \dfrac{1}{x} |xR_{V,0}(x)|^{q'}\,dx<+\infty,\ q'\geq 1,\quad}\\
\\
{ \gamma_0^2< k^2-1/4\quad}
\end{array}
\right.
\end{equation}
and
\begin{equation} \label{ipotesi-potenziale-zero-resto2}
 {\mbox{if $\mu_a\neq 0$}}: \quad \left\{\begin{array}{l}
\displaystyle{x^{\alpha_0} R_{V,0}=o(1),\quad
 x^{\alpha_0+1} R'_{V,0}=o(1),\ x\to 0^+,\quad }\\
\\
\displaystyle{\int_0^1 \dfrac{1}{x^{\alpha_0+1}} |x^{\alpha_0+1}R'_{V,0}(x)|^{q''}\,dx<+\infty,\ q''\geq 1,\quad}\\
\\
\gamma_0\neq 0.\quad
\end{array}
\right.
\end{equation}
Under these conditions, assumption $({\cal P}_1)$ is satisfied with $\mu^\pm=\pm 1$; indeed, we obviously have
\[
\lim_{x\to +\infty} P_{V,k,\mu_a} (x)=\lim_{x\to +\infty} \left(\begin{array}{cc}
 	      -1+V(x)&{\displaystyle{-\dfrac{k}{x}-\mu_a V'(x)}}\\
&\\
			         {\displaystyle{-\dfrac{k}{x}}-\mu_a V'(x)}&1+V(x)
				\end{array}
			\right)=\left(\begin{array}{cc}
 	      -1&0\\
&\\
			        0&1
				\end{array}
			\right).
\]
Moreover, the matrix $R_\infty$ in \eqref{iporestoinfinito} is given by
\[
R_\infty (x)=\left(\begin{array}{cc}
 	      V(x)&{\displaystyle{-\dfrac{k}{x}-\mu_a V'(x)}}\\
&\\
			         {\displaystyle{-\dfrac{k}{x}}-\mu_a V'(x)}&V(x)
				\end{array}
			\right),\quad \forall \  x>0;
\]
from \eqref{ipotesi-potenziale-infinito} we deduce that $V, R_{V,\infty} \in L^q(1,+\infty)$, for every $q>1/\alpha_\infty$, and $R'_{V,\infty}\in L^p(1,+\infty)$, for every $p>1/(\alpha_\infty+1)$, while we plainly have $k/x\in L^s(1,+\infty)$, for every $s>1$. By observing that all the functions $V, R_{V,\infty},  R'_{V,\infty}, k/x$ go to zero at infinity, we conclude that \eqref{iporestoinfinito} is satisfied with $q_\infty=1/(\alpha_\infty+1)$.
\vs{6}
\ni
As far as $({\cal P}_2)$ and $({\cal P}_3)$ are concerned, a crucial role is played by the constants $\gamma_0$ and $\mu_a$, as it is evident from the assumptions on $V$. Indeed, let us first discuss the case $\mu_a=0$; in this situation, taking $\alpha_0=1$ and $\beta=\alpha_0=1$ in \eqref{limitiPzero}, we have
\[
\begin{array}{l}
\displaystyle{P^*=\lim_{x\to 0^+} x  P_{V,k,0}(x)=}\\
\\
\displaystyle{=\lim_{x\to 0^+} \left(\begin{array}{cc}
 	      -x+xV(x)&{\displaystyle{-k}}\\
&\\
			         {\displaystyle{-k}}&x+xV(x)
				\end{array}
			\right)= \left(\begin{array}{cc}
 	      \gamma_0&{\displaystyle{-k}}\\
&\\
			         {\displaystyle{-k}}&\gamma_0
				\end{array}
			\right)}
\end{array}
\]
and
\[
R_0(x)=xP(x)-P^*=\left(\begin{array}{cc}
 	     -x+xV(x)-\gamma_0&0\\
&\\
			         0&x+xV(x)-\gamma_0
				\end{array}
			\right).
\]
Hence, from \eqref{ipotesi-potenziale-zero} and \eqref{ipotesi-potenziale-zero-resto1} we infer
\[
\int_0^1 \dfrac{1}{x} \ ||R_0(x)||^{q'}\,dx \leq 2^{q'-1}\int_0^1 (x^{q'-1}+x^{q'-1} |R_{V,0}(x)|^{q'})\,dx<+\infty,
\]
concluding that \eqref{iporestozero} holds true with $q_0=q'$. Moreover, the last relation in \eqref{ipotesi-potenziale-zero-resto1} guarantees that \eqref{determinante-beta-1} is fulfilled.
\vs{6}
\ni
Suppose now $\mu_a\neq 0$; taking $\beta=\alpha_0+1$ in \eqref{limitiPzero}, we have
\[
\begin{array}{l}
\displaystyle{P^*=\lim_{x\to 0^+} x^{\alpha_0+1}  P_{V,k,\mu_a}(x)=}\\
\\
\displaystyle{=\lim_{x\to 0^+} \left(\begin{array}{cc}
 	      -x^{\alpha_0+1}+x^{\alpha_0+1}V(x)&{\displaystyle{-kx^{\alpha_0}-\mu_a x^{\alpha_0+1}V'(x)}}\\
&\\
			         {\displaystyle{-kx^{\alpha_0}-\mu_a x^{\alpha_0+1}V'(x)}}&x^{\alpha_0+1}+x^{\alpha_0+1}V(x)
				\end{array}
			\right)= \left(\begin{array}{cc}
 	      0&{\displaystyle{\mu_a \alpha_0 \gamma_0}}\\
&\\
			         {\displaystyle{\mu_a \alpha_0 \gamma_0}}&0
				\end{array}
			\right)}
\end{array}
\]
and
\[
R_0(x)=x^{\alpha_0+1}  P(x)-P^*=\left(\begin{array}{cc}
 	      -x^{\alpha_0+1}+x^{\alpha_0+1}V(x)&-kx^{\alpha_0}-\mu_a x^{\alpha_0+1}V'(x)-\mu_a \alpha_0 \gamma_0\\
&\\
			         -kx^{\alpha_0}-\mu_a x^{\alpha_0+1}V'(x)-\mu_a \alpha_0 \gamma_0& x^{\alpha_0+1}+x^{\alpha_0+1}V(x)
				\end{array}
			\right).
\]
Now, let $q_0>\max (q'',\alpha_0)$; from \eqref{ipotesi-potenziale-zero} and \eqref{ipotesi-potenziale-zero-resto2} we infer
\[
\begin{array}{l}
\displaystyle{\int_0^1 \dfrac{1}{x^{\alpha_0+1}}  |x^{\alpha_0+1}|^{q_0}\,dx =\int_0^1 \dfrac{1}{x^{(\alpha_0+1)(1-q_0)}}\,dx<+\infty}\\
\\
\displaystyle{\int_0^1 \dfrac{1}{x^{\alpha_0+1}}  |x^{\alpha_0+1}V(x)|^{q_0}\,dx \leq 2^{q_0-1} \int_0^1 \dfrac{x^{q_0}}{x^{\alpha_0+1}} \left(\gamma_0^{q_0}+|x^{\alpha_0}R_{V,0}(x)|^{q_0}\right)\,dx<+\infty}\\
\\
\displaystyle{\int_0^1 \dfrac{1}{x^{\alpha_0+1}}  |k x^{\alpha_0}|^{q_0}\,dx =|k|^{q_0} \ \int_0^1 \dfrac{1}{x^{1+(1-q_0)\alpha_0}}\,dx<+\infty}\\
\\
\displaystyle{\int_0^1 \dfrac{1}{x^{\alpha_0+1}}  |-\mu_a x^{\alpha_0+1}V'(x)-\mu_a\alpha_0 \gamma_0|^{q_0}\,dx =}\\
\\
\displaystyle{\qquad \qquad \qquad \qquad =|\mu_a|^{q_0} \int_0^1 \dfrac{1}{x^{\alpha_0+1}} |x^{\alpha_0+1}R'_{V,0}(x)|^{q''}  |x^{\alpha_0+1}R'_{V,0}(x)|^{q_0-q''} \,dx<+\infty,}
\end{array}
\]
concluding again that \eqref{iporestozero} holds true. Moreover, the last relation in \eqref{ipotesi-potenziale-zero-resto2} guarantees that \eqref{determinante-beta} is fulfilled.
\vs{12}
\ni
The fact that $P\in {\cal P}_{\mu} $, and in particular the fact that $P$ satisfies $({\cal P}_3)$, is related to the spectral properties of the operator $\tau: z\to Jz'+P(x)z$; indeed, as we will see at the beginning of \ref{sez-autovalori}, condition $({\cal P}_3)$ implies that the operator $\tau$ is in the limit point case at $x=0$. As a consequence, it admits a unique self-adjoint realization (cfr. \eqref{def-A0}); in the particular case of the operator associated to \eqref{Dirac-magn}, with the Coulomb potential
\begin{equation} \label{pot-Coulomb}
V(x)=\dfrac{\gamma}{x},\quad \forall \ x>0,\ \gamma<0,
\end{equation}
condition $({\cal P}_3)$ is satisfied for a larger range of values of $\gamma$ when $\mu_a\neq 0$. This means that the presence of an anomalous magnetic moment has a regularizing effect on the radial Dirac operator (cf. also \cite[Sect. 5.3.2]{Th-book-92}.

\end{remark}
\vs{12}
\ni
\subsection{Asymptotic estimates} \label{stimeasintotiche}
\vs{12}
\ni
In this subsection we describe the behaviour of the solutions of \eqref{lineare1} when $x\to +\infty$ or $x\to 0^+$; this will be the consequence of some general results on the asymptotic properties of solutions of linear equations (see e.g. \cite{Ea-book-89}). As a first step, let us consider a system of the form
\begin{equation} \label{lineare-Levinson}
u'=C(\l)u+U(x,\l)u,\quad x\geq 1,\l \in \Lambda
\end{equation}
where $C(\l)$ and $U(x,\l)$ are $2\times2$ matrix, for every $\l\in \Lambda$ and $x\geq 1$. We have the following result:
\vs{12}
\ni
\begin{proposition} \label{Levinson-1} (\cite[Th. 1.5.2, Th. 1.8.1, Th. 1.8.2]{Ea-book-89}) Let us suppose that for every $\l\in \Lambda$ the matrix $C(\l)$ has two real eigenvalues $\sigma^-_\l< 0 < \sigma^+_\l$ and let $u^-_\l$, $u^+_\l$ be the eigenvectors associated to $\sigma^-_\l$ and $\sigma^+_\l$, respectively. Moreover, let us assume that
\begin{equation} \label{resto-infinitesimo}
\lim_{x\to +\infty} U(x,\l)=0,\quad \forall \ \l \in \Lambda,
\end{equation}
and that there exists $q\geq 1$ such that
\begin{equation} \label{iporesto-Levinson}
\int_1^{+\infty} ||U(x,\l)||^q\,dx<+\infty,\quad \forall \ \l \in \Lambda.
\end{equation}
Then, for every $\l \in \Lambda$ system \eqref{lineare-Levinson} has two linearly independent solutions $u_{1,\l}$ and $u_{2,\l}$ satisfying
\begin{equation} \label{stimeu}
\begin{array}{l}
\displaystyle{u_{1,\l}(x)=(u^-_\l+o(1))e^{\sigma^-_\l (x-1)+\int_1^x g_{1,\l}(t)\,dt},\ x\to +\infty}\\
\\
\displaystyle{u_{2,\l}(x)=(u^+_\l+o(1))e^{\sigma^+_\l (x-1)+\int_1^x g_{2,\l}(t)\,dt},\ x\to +\infty,}
\end{array}
\end{equation}
where, for $i=1, 2$, we have
\begin{equation} \label{espressione-resti}
\begin{array}{ll}
                       g_{i,\l}=0&{\mbox{if $q=1$}}\\
                        &\\
                        g_{i,\l}\in L^q (1,+\infty)&{\mbox{if $q>1$.}}
                       \end{array}
\end{equation}
\end{proposition}
\vs{8}
\ni
{\sc Proof.} Let us note that when $q=1$ the result follows from \cite[Th. 1.8.1]{Ea-book-89}. Therefore, assume that $q>1$; from \cite[Th. 1.5.2, Th. 1.8.2]{Ea-book-89} we immediately deduce that \eqref{stimeu} is satisfied with some functions $g_{i,\l}$, $i=1, 2$, $\l \in \Lambda$, such that
\begin{equation} \label{espressione-resti-1}
g_{i,\l}=\left\{\begin{array}{ll}
                       0&{\mbox{if $q=1$}}\\
                        &\\
                        \displaystyle{ \sum_{m=1}^M g_{i,m,\l}}&{\mbox{if $q>1$,}}
                       \end{array}
              \right.
\end{equation}
with $M$ such that $2^{M-1}<q\leq 2^M$ and
\begin{equation} \label{sommabilita-resti}
\begin{array}{l}
g_{i,m,\l}\in L^{q/2^{m-1}}(1,+\infty),\quad \forall m=1,\ldots , M.
\end{array}
\end{equation}
Now, assumption \eqref{resto-infinitesimo} implies that
\[
\lim_{x\to +\infty} g_{i,m,\l} (x)=0,\quad \forall m=1,\ldots , M,\ \l \in \Lambda,\ i=1, 2
\]
(see also formula (1.5.27) in \cite{Ea-book-89}). Hence, for every $i=1, 2$, $\l \in \Lambda$ and $m=1,\ldots , M$ we have
\[
g_{i,m,\l}\in  L^{q/2^{m-1}}(1,+\infty)\quad \Rightarrow \quad g_{i,m,\l} \in L^q(1,+\infty).
\]
This implies that $q_{i,\l}\in L^q (1,+\infty)$, for every $i=1, 2$ and $\l \in \Lambda$.
\qed
\vs{12}
\ni
Now, let us observe that
\begin{equation} \label{stima-primitiva}
f\in L^p(1,+\infty),\ p>1 \quad \Rightarrow \quad \left|\int_1^x f(t)\,dt\right|\leq \| f \|_{ L^{p} }  (x-1)^{1/p'},\quad \forall \ x\geq 1,
\end{equation}
where $p'$ is the conjugate exponent of $p$; noting that in this case $1/p'<1$, from Proposition \ref{Levinson-1} we obtain the following result:
\vs{12}
\ni
\begin{proposition} \label{Levinson-2} Under the assumptions of Proposition \ref{Levinson-1}, for every $\l \in \Lambda$ we have
\begin{equation} \label{limiteu1}
\lim_{x\to +\infty} u_{1,\l}(x)=0
\end{equation}
and
\begin{equation} \label{limiteu2}
\lim_{x\to +\infty} \|u_{2,\l}(x)\|=+\infty. 
\end{equation}
Moreover, if $q>1$ in \eqref{iporesto-Levinson}, then there exists $x_1>1$ such that
\begin{equation} \label{stima-esponenziale}
\sigma^-_\l (x-1)+\int_1^x g_{1,\l}(t)\,dt\leq \dfrac{\sigma^-_\l}{2} (x-1),\quad \forall \ x\geq x_1.
\end{equation}
\end{proposition}
\vs{12}
\ni
Using Proposition \ref{Levinson-1} and Proposition \ref{Levinson-2} we are able to prove some asymptotic results on the solutions of \eqref{lineare1} when $x\to +\infty$ or $x\to 0^+$. We start with the study of \eqref{lineare1} when $x\to +\infty$ (cf. also \cite{CaDa-10}); assume then $x\geq 1$.
\vs{4}
\ni
Let us first observe that \eqref{lineare1} can be written as
\begin{equation} \label{sistemaintermedio}
\begin{array}{l}
                         z'=B_\l z+Q(x)z,
\end{array}
\end{equation}
where
\[
B_\l=J^{-1} (\l {\rm Id}- P_\infty),\quad Q(x)=J^{-1} (P_\infty-P(x)),\quad \forall \ x>0.
\]
This form corresponds to \eqref{lineare-Levinson} with $C(\l)=B_\l$ and $U(x,\l)=Q(x)$, for every $x\geq 1$, $\l \in \Lambda$; note that assumptions \eqref{limitiPinfinito} and \eqref{iporestoinfinito} imply that \eqref{resto-infinitesimo} and \eqref{iporesto-Levinson}, with $q=q_\infty$, hold true. Moreover, if $\l \in \Lambda$, setting $\Delta_{\lambda}=(\mu^+-\lambda)(\lambda-\mu^-)$, then $B_\l$ has the real eigenvalues $\pm \sqrt{\Delta_\l}$; in this situation we denote by $b_{1,\l} = (\lambda- \mu^{+}  , \sqrt{ \Delta_\l } ) $ and $b_{2,\l}  = (\mu^{+} - \lambda , \sqrt{ \Delta_\l } ) $ the eigenvectors of $B_\l$ associated to the eigenvalues $-\sqrt{\Delta_\l}$ and $\sqrt{\Delta_\l}$, respectively.
\vs{4}
\ni
Therefore from Proposition \ref{Levinson-1} and Proposition \ref{Levinson-2} we deduce the following results:
\vs{12}
\ni
\begin{proposition} \label{stimeinf2} For every $\l \in \Lambda$ system \eqref{lineare1} has two linearly independent solutions $z_{1,\l}$ and $z_{2,\l}$ satisfying
\begin{equation} \label{stimez}
\begin{array}{l}
\displaystyle{z_{1,\l}(x)=(b_{1,\l}+o(1))e^{-\sqrt{\Delta_\l}(x-1)+\int_1^x g_{1}(t)\,dt},\ x\to +\infty}\\
\\
\displaystyle{z_{2,\l}(x)=(b_{2,\l}+o(1))e^{\sqrt{\Delta_\l}(x-1)+\int_1^x g_{2}(t)\,dt},\ x\to +\infty,}
\end{array}
\end{equation}
where, for $i=1, 2$, we have
\begin{equation} \label{espressione-resti-infinito}
\begin{array}{ll}
                       g_{i}=0&{\mbox{if $q_\infty=1$}}\\
                        &\\
                        g_{i}\in L^{q_\infty} (1,+\infty)&{\mbox{if $q_\infty>1$.}}
                       \end{array}
\end{equation}
\end{proposition}
\vs{12}
\ni
\begin{lemma} \label{limitiuv1} Assume that $\l\in \Lambda$ and let $z_{1,\l}$ and $z_{2,\l}$ be the solutions of \eqref{lineare1} given in Proposition \ref{stimeinf2}. Then
\begin{equation} \label{limitez1}
\lim_{x\to +\infty} z_{1,\l}(x)=0
\end{equation}
and
\begin{equation} \label{limitez2}
\lim_{x\to +\infty} |(z_{2,\l})_1(x)|=\lim_{x\to +\infty} |(z_{2,\l})_2(x)|=+\infty.
\end{equation}
Moreover, $z_{1,\l}\in H^1(1,+\infty)$.
\end{lemma}
\vs{8}
\ni
{\sc Proof.} The relations \eqref{limitez1} and \eqref{limitez2} immediately follow from \eqref{limiteu1} and \eqref{limiteu2}.
In particular \eqref{limitez2} comes from the fact that neither component of $ b_{ 2, \lambda } $ vanishes.
\vs{4}
\ni
Moreover, from \eqref{stima-esponenziale} we deduce that there exists $K_{1,\l}>0$ such that
\[
||z_{1,\l}(x)||\leq K_{1,\l} e^{-\sqrt{\Delta_\l}(x-1)},\quad \forall \ x\geq x_1;
\]
this implies that $z_{1,\l}\in L^2(1,+\infty)$. Now, from the differential equation we deduce that
\[
Jz'_{1,\l}(x)=\l z_{1,\l}(x)-P(x)z_{1,\l}(x),\quad \forall \ x\geq 1;
\]
since $P\in L^\infty(1,+\infty)$, we infer that $Jz'_{1,\l}\in L^2(1,+\infty)$ and then $z_{1,\l}\in H^1(1,+\infty)$.
\qed
\vs{12}
\ni
Arguing as in the proof of \cite[Lemma 2.3]{CaDa-10}, we obtain the following result:
\vs{12}
\ni
\begin{lemma} \label{limitiuv} Assume that $\l\in \Lambda$ and let $z=(u,v)$ be a nontrivial solution of \eqref{lineare1}. Then either
\begin{equation} \label{limiteuv1}
\lim_{x\to +\infty} u(x)=\lim_{x\to +\infty} v(x)=0
\end{equation}
or
\begin{equation} \label{limiteuv2}
\lim_{x\to +\infty} |u(x)|=\lim_{x\to +\infty} |v(x)|=+\infty.
\end{equation}
Moreover, $z\in H^1(1,+\infty)$ if and only if \eqref{limiteuv1} holds true and there exists $\gamma>0$ such that $z=\gamma z_{1,\l}$, where $z_{1,\l}$ is given in Proposition \ref{stimeinf2}.
\end{lemma}
\vs{12}
\ni
Now, let us study the behaviour of the solutions of \eqref{lineare1} when $x\to 0^+$; assume then that $x\in (0,1)$. For every $\beta\geq 1$ let us consider an invertible function $\phi_\beta\in C^1((1,+\infty),(0,1))$ such that
\begin{equation} \label{ipotesi-phi-beta}
\lim_{t\to +\infty} \phi_\beta(t)=0 \quad \text{and} \quad \lim_{ t \to 1^{+} } \phi_{\beta}( t ) = 1 .
\end{equation}
The change of variable $x=\phi_\beta (t)$ transforms \eqref{lineare1} into
\begin{equation} \label{sistematrasformato-generale}
w'=-J^{-1} P(\phi_\beta (t)) \phi'_\beta (t)w+\l J^{-1}  \phi'_\beta (t) w,
\end{equation}
where $w(t)=z(\phi_\beta (t))$, for every $t\geq 1$. With a suitable choice of $\phi_\beta$ system \eqref{sistematrasformato-generale} can be reduced to a system of the form \eqref{lineare-Levinson}:
\vs{12}
\ni
\begin{lemma} \label{scelta-phi-beta-1} Assume $\beta=1$ in $({\cal P}_2)$ and let
\[
\phi_\beta (t)=e^{1-t} ,\quad \forall \ t\geq 1.
\]
Then \eqref{sistematrasformato-generale} reduces to a system of the form \eqref{lineare-Levinson} with
\begin{equation} \label{espressioni-matrici-1}
C=C(\l)=J^{-1} P^*,\quad U(t,\l)=J^{-1} (R_0( e^{1-t})-\l  e^{1-t} {\rm Id}),
\end{equation}
for every $t\geq 1$ and $\l \in \Lambda$.
\end{lemma}
\vs{12}
\ni
\begin{lemma} \label{scelta-phi-beta-gen} Assume $\beta>1$ in $({\cal P}_2)$  and let
\[
\phi_\beta (t)=t^{-1/(\beta -1)},\quad \forall \ t\geq 1.
\]
Then \eqref{sistematrasformato-generale} reduces to a system of the form \eqref{lineare-Levinson} with
\begin{equation} \label{espressioni-matrici-gen}
C=C(\l)=\dfrac{1}{\beta -1} J^{-1} P^*,\quad U(t,\l)=\dfrac{1}{\beta -1}  J^{-1} (R_0(t^{-1/(\beta -1)})-\l t^{-\beta/(\beta -1)}{\rm Id}),
\end{equation}
for every $t\geq 1$ and $\l \in \Lambda$.
\end{lemma}
\vs{12}
\ni
The proofs of Lemma \ref{scelta-phi-beta-1} and Lemma \ref{scelta-phi-beta-gen} are straightforward and therefore they are omitted.
\vs{12}
\ni
Now, set $\Delta^* = -{\rm det }\  P^*$ and observe that the matrix $C$ given in \eqref{espressioni-matrici-1} or \eqref{espressioni-matrici-gen} has the eigenvalues $\sigma^\pm = \pm \sqrt{\Delta^*}$ if $\beta=1$ and $\sigma^\pm = \pm \sqrt{\Delta^*} / (\beta -1) $ if $\beta>1$; in what follows, we will denote by $w_1^*$ and $w_2^*$ the eigenvectors of $C$ associated to $\sigma^\pm$.
\vs{4}
\ni
Moreover, from \eqref{limitiPzero} and the definition of $R_0$ we deduce that the function $U$ given in \eqref{espressioni-matrici-1} or \eqref{espressioni-matrici-gen} satisfies \eqref{resto-infinitesimo}. Finally, let us note that \eqref{iporestozero} implies that \eqref{iporesto-Levinson} is satisfied with $q=q_0$; indeed, when $\beta=1$ we have
\[
\begin{array}{l}
{\displaystyle{\int_1^{+\infty} ||U(t,\l)||^{q_0}\,dt \leq \left\{ \left[ \int_1^{+\infty} ||R_0(e^{1-t})||^{q_0}\,dt \right]^{\frac{1}{q_{0}}} + \l \left[ \int_1^{+\infty} e^{q_0(1-t)}\,dt \right]^{\frac{1}{q_{0}}} \right\}^{q_{0}} =}} \\
\\
{\displaystyle{= \left\{ \left[ \int_0^{1} \dfrac{1}{x}\ ||R_0(x)||^{q_0}\,dx \right]^{\frac{1}{q_{0}}} + \frac{\lambda}{q_{0}^{1/q_{0}}} \right\}^{q_{0}}<+\infty.}}
\end{array}
\]
On the other hand, if $\beta>1$ we deduce that
\[
\begin{array}{l}
{\displaystyle{\int_1^{+\infty} ||U(t,\l)||^{q_0}\,dt\leq \dfrac{1}{(\beta -1)^{q_{0}}} \left\{ \left[ \int_1^{+\infty} ||R_0(t^{-1/(\beta -1)})||^{q_0}\,dt \right]^{\frac{1}{q_{0}}} + \l \left[ \int_1^{+\infty} t^{-\beta q_0/(\beta -1)}\,dt \right]^{\frac{1}{q_{0}}} \right\}^{q_{0}} = }}\\
\\
{\displaystyle{= \dfrac{1}{\beta -1} \left\{ \left[ (\beta -1) \int_0^{1} \dfrac{1}{x^\beta}\ ||R_0(x)||^{q_0}\,dx \right]^{\frac{1}{q_{0}}} + \l \left(\frac{\beta-1}{q_{0}}\right)^{\frac{1}{q_{0}}} \right\}^{q_{0}} <+\infty.}}
\end{array}
\]
Therefore, we can apply Proposition \ref{Levinson-1} and Proposition \ref{Levinson-2} to \eqref{sistematrasformato-generale}, with $\phi_\beta$ as above, and obtain the following results:
\vs{12}
\ni
\begin{proposition} \label{stimezero2-0} For every $\l \in \Lambda$ system \eqref{sistematrasformato-generale}, with $\phi_\beta$ as in Lemma \ref{scelta-phi-beta-1} or Lemma \ref{scelta-phi-beta-gen}, has two linearly independent solutions $w_{1,\l}$ and $w_{2,\l}$ satisfying
\begin{equation} \label{stimew-1}
\begin{array}{l}
\displaystyle{w_{1,\l}(t)=(w_1^*+o(1))e^{-\sqrt{\Delta^*}(t-1)+\int_1^t g_{1,\l}(s)\,ds},\ t\to +\infty}\\
\\
\displaystyle{w_{2,\l}(t)=(w_2^*+o(1))e^{\sqrt{\Delta^*}(t-1)+\int_1^t g_{2,\l}(s)\,ds},\ t\to +\infty,}
\end{array}
\end{equation}
if $\beta =1$ and
\begin{equation} \label{stimew-gen}
\begin{array}{l}
\displaystyle{w_{1,\l}(t)=(w_1^*+o(1))e^{-\frac{\sqrt{\Delta^*}}{\beta -1}  (t-1)+\int_1^t g_{1,\l}(s)\,ds},\ t\to +\infty}\\
\\
\displaystyle{w_{2,\l}(t)=(w_2^*+o(1))e^{\frac{\sqrt{\Delta^*}}{\beta -1}  (t-1)+\int_1^t g_{2,\l}(s)\,ds},\ t\to +\infty,}
\end{array}
\end{equation}
if $\beta >1$, where, for $i=1, 2$, we have
\begin{equation} \label{espressione-resti-zero}
\begin{array}{ll}
                       g_{i,\l}=0&{\mbox{if $q_0=1$}}\\
                        &\\
                        g_{i,\l}\in L^{q_0} (1,+\infty)&{\mbox{if $q_0>1$.}}
                       \end{array}
\end{equation}
\end{proposition}
\vs{12}
\ni
\vs{12}
\ni
\begin{lemma} \label{limitiuv2} Assume that $\l\in \Lambda$ and let $w_{1,\l}$ and $w_{2,\l}$ be the solutions of \eqref{sistematrasformato-generale} given in Proposition \ref{stimezero2-0}. Then
\begin{equation} \label{limitew1}
\lim_{t\to +\infty} w_{1,\l}(t)=0
\end{equation}
and
\begin{equation} \label{limitew2}
\lim_{t\to +\infty} \|w_{2,\l}(t)\| =+\infty.
\end{equation}
Moreover, the solution $w_{1,\l}$ satisfies
\begin{equation} \label{lemma-w-1}
\displaystyle{\int_{1}^{+\infty} ||w_{1,\l}(t)||^2\ e^{t}\,dt<+\infty,\quad {\mbox{if $\beta=1$}}}
\end{equation}
and
\begin{equation} \label{lemma-w-gen}
\displaystyle{\int_{1}^{+\infty} ||w_{1,\l}(t)||^2\ t^{\beta/(\beta-1)}\,dt<+\infty,\quad {\mbox{if $\beta>1.$}}}
\end{equation}
\end{lemma}
\vs{8}
\ni
{\sc Proof.} Let us note that \eqref{limitew1} and \eqref{limitew2} immediately follow from \eqref{limiteu1} and \eqref{limiteu2}.
\vs{6}
\ni
As far as \eqref{lemma-w-1} is concerned, from \eqref{stimew-1} we deduce that there exists $K_{1,\l}>0$ such that
\begin{equation} \label{2nuovo}
||w_{1,\l}(t)||^2\ e^{t} \sim K_{1,\l} e^{-2\sqrt{\Delta^*}(t-1)+2\int_1^t g_{1,\l}(s)\,ds}\ e^{t},\quad t\to +\infty;
\end{equation}
now, let us observe that $1-2\sqrt{\Delta^*}>0$, since \eqref{determinante-beta-1} holds. Hence, using again \eqref{stima-primitiva} we infer that there exists $t_1>1$ such that
\begin{equation} \label{1nuovo}
e^{(1-2\sqrt{\Delta^*)}t+2\int_1^t g_{1,\l}(s)\,ds}\leq e^{(1-2\sqrt{\Delta^*})t/2},\quad \forall \ t\geq t_1,
\end{equation}
is satisfied. Conditions \eqref{2nuovo} and \eqref{1nuovo} imply \eqref{lemma-w-1}.
\vs{6}
\ni
Finally, when $\beta>1$ from \eqref{stimew-gen} we deduce that there exists $M_{1,\l}>0$ such that
\begin{equation} \label{3nuovo}
||w_{1,\l}(t)||^2\ t^{\beta/(\beta-1)} \sim M_{1,\l} e^{-2 \frac{\sqrt{\Delta^*}}{\beta -1}(t-1)+2\int_1^t g_{1,\l}(s)\,ds}\  t^{\beta/(\beta-1)},\quad t\to +\infty;
\end{equation}
moreover, from \eqref{stima-esponenziale} we infer that there exists $t_2>1$ such that
\begin{equation} \label{4nuovo}
e^{-2 \frac{\sqrt{\Delta^*}}{\beta -1}  (t-1)+2\int_1^t g_{1,\l}(s)\,ds}\leq e^{-  \frac{\sqrt{\Delta^*}}{\beta -1}  (t-1)},\quad \forall \ t\geq t_2,
\end{equation}
is satisfied. Conditions \eqref{3nuovo} and \eqref{4nuovo} imply \eqref{lemma-w-gen}.
\qed
\vs{12}
\ni
The next result is a consequence of Proposition \ref{stimezero2-0} and Lemma \ref{limitiuv2}.
\vs{12}
\ni
\begin{proposition} \label{stimezero2} For every $\l \in \Lambda$ system \eqref{lineare1} has two linearly independent solutions $\zeta_{1,\l}$ and $\zeta_{2,\l}$ satisfying
\begin{equation} \label{stimezeta-1}
\begin{array}{l}
\displaystyle{\zeta_{1,\l}(x)=(w_1^*+o(1))\ x^{\sqrt{\Delta^*}}e^{\int_1^{1-\log x} g_{1,\l} (s)\,ds},\ x\to 0^+}\\
\\
\displaystyle{\zeta_{2,\l}(x)=(w_2^*+o(1))\  x^{-\sqrt{\Delta^*}}  e^{\int_1^{1-\log x} g_{2,\l} (s)\,ds},\ x\to 0^+,}
\end{array}
\end{equation}
if $\beta=1$ and
\begin{equation} \label{stimezeta-gen}
\begin{array}{l}
\displaystyle{\zeta_{1,\l}(x)=(w_1^*+o(1))\ e^{- \frac{\sqrt{\Delta^*}}{\beta-1}  x^{1-\beta} +\int_1^{-\log x} g_{1,\l} (s)\,ds},\ x\to 0^+}\\
\\
\displaystyle{\zeta_{2,\l}(x)=(w_2^*+o(1))\ e^{  \frac{\sqrt{\Delta^*}}{\beta-1}  x^{1-\beta} \int_1^{-\log x} g_{2,\l} (s)\,ds},\ x\to 0^+,}
\end{array}
\end{equation}
if $\beta >1$, where, for $i=1, 2$, we have
\begin{equation} \label{espressione-resti-zero-2}
\begin{array}{ll}
                       g_{i,\l}=0&{\mbox{if $q_0=1$}}\\
                        &\\
                        g_{i,\l}\in L^{q_0} (1,+\infty)&{\mbox{if $q_0>1$.}}
                       \end{array}
\end{equation}
\end{proposition}
\vs{12}
\ni
\begin{lemma} \label{limitiuvzero1} Assume that $\l\in \Lambda$ and let $\zeta_{1,\l}$ and $\zeta_{2,\l}$ be the solutions of \eqref{lineare1} given in Proposition \ref{stimezero2}. Then
\begin{equation} \label{limitezeta1}
\lim_{x\to 0^+} \zeta_{1,\l}(x)=0
\end{equation}
and
\begin{equation} \label{limitezet2}
\lim_{x\to 0^+} \|\zeta_{2,\l}(x)\| =+\infty.
\end{equation}
Moreover, $\zeta_{1,\l}\in H^1(0,1)$.
\end{lemma}
\vs{8}
\ni
{\sc Proof.} The relations \eqref{limitezeta1} and \eqref{limitezet2} immediately follow from \eqref{limitew1} and \eqref{limitew2}.
\vs{6}
\ni
Now, assume that $\beta=1$; let us observe that we have
\begin{equation} \label{zetaprimo1l2}
\int_0^1 \dfrac{||\zeta_{1,\l}(x)||^2}{x^2}\,dx=\int_{1}^{+\infty} ||\zeta_{1,\l} (e^{1-t})||^2 e^{t-1}\,dt=\int_{1}^{+\infty} ||w_{1,\l} (t)||^2 e^{t-1}\,dt  <+\infty,
\end{equation}
by \eqref{lemma-w-1}. This condition obviously implies that
\begin{equation} \label{zeta1l2}
\int_0^1 ||\zeta_{1,\l} (x)||^2\,dx\leq \int_0^1 \dfrac{||\zeta_{1,\l}(x)||^2}{x^2}\,dx<+\infty
\end{equation}
and so $\zeta_{1,\l}\in L^2(0,1)$; in order to prove that $\zeta'_{1,\l}\in L^2(0,1)$, let us note that $\zeta_{1,\l}$ satisfies
\[
J\zeta'_{1,\l}(x)=\l \zeta_{1,\l}(x)-P(x)\zeta_{1,\l}(x),\quad \forall \ x\in (0,1).
\]
From $({\cal P}_2)$ we deduce that there exists $K>0$ such that
\[
||P(x)||\leq \dfrac{K}{x},\quad \forall \ x\in (0,1);
\]
therefore, \eqref{zetaprimo1l2} implies that $P\zeta_{1,\l}\in L^2(0,1)$ and then also $\zeta'_{1,\l}\in L^2(0,1)$.
\vs{6}
\ni
When $\beta>1$ by \eqref{lemma-w-gen} we obtain
\begin{equation} \label{zetaprimo1l2-gen}
\int_0^1 \dfrac{||\zeta_{1,\l}(x)||^2}{x^{2\beta}}\,dx=\dfrac{1}{\beta -1} \ \int_{1}^{+\infty} \dfrac{||\zeta_{1,\l} (t^{-1/(\beta -1)})||^2}{t^{-2\beta /(\beta -1)}} \ t^{-\beta/(\beta -1)}\,dt=\int_{1}^{+\infty} ||w_{1,\l} (t)||^2\ t^{\beta/(\beta -1)} \,dt<+\infty.
\end{equation}
This condition implies that
\begin{equation} \label{zeta1l2-gen}
\int_0^1 ||\zeta_{1,\l} (x)||^2\,dx\leq \int_0^1 \dfrac{||\zeta_{1,\l}(x)||^2}{x^{2\beta}}\,dx<+\infty
\end{equation}
and so $\zeta_{1,\l}\in L^2(0,1)$; arguing as above, in order to prove that $\zeta'_{1,\l}\in L^2(0,1)$, let us note that $\zeta_{1,\l}$ satisfies
\[
J\zeta'_{1,\l}(x)=\l \zeta_{1,\l}(x)-P(x)\zeta_{1,\l}(x),\quad \forall \ x\in (0,1).
\]
From $({\cal P}_2)$ we deduce that there exists $M>0$ such that
\[
||P(x)||\leq \dfrac{M}{x^{\beta}},\quad \forall \ x\in (0,1);
\]
therefore, \eqref{zetaprimo1l2-gen} implies that $P\zeta_{1,\l}\in L^2(0,1)$ and then also $\zeta'_{1,\l}\in L^2(0,1)$.
\qed
\vs{12}
\ni
\begin{lemma} \label{limitiuvzero} Assume that $\l\in \Lambda$ and let $z=(u,v)$ be a nontrivial solution of \eqref{lineare1}. Then either
\begin{equation} \label{limiteuvzero1}
\lim_{x\to 0^+} u(x)=\lim_{x\to 0^+} v(x)=0
\end{equation}
or
\begin{equation} \label{limiteuvzero2}
\lim_{x\to 0^+} \|z(x)\|  =+\infty.
\end{equation}
Moreover, $z\in H^1(0,1)$ if and only if \eqref{limiteuvzero1} holds true and there exists $\xi \in \R  $ such that $z=\xi \zeta_{1,\l}$, where $\zeta_{1,\l}$ is given in Proposition \ref{stimezero2}.
\end{lemma}
\vs{12}
\ni
\begin{remark} \label{spazio-soluzioni} Let us denote by $Z$ the set of solutions of \eqref{lineare1}. From Proposition \ref{stimeinf2} and Proposition \ref{stimezero2} we deduce that
\[
Z={\rm span } \  \{z_{1,\l}, z_{2,\l}\}={\rm span } \  \{\zeta_{1,\l}, \zeta_{2,\l}\}.
\]
As far as nontrivial solutions $z\in L^2(0,+\infty)$ are concerned, let us observe that  Lemma \ref{limitiuv} and Lemma \ref{limitiuvzero} prove that
\[
z\in L^2(1,+\infty) \quad \Longleftrightarrow \quad z\in {\rm span } \  \{z_{1,\l}\}:=Z_\infty
\]
and
\[
z\in L^2(0,1) \quad \Longleftrightarrow \quad z\in {\rm span } \  \{\zeta_{1,\l}\}:=Z_0.
\]
As a consequence, $z\in L^2(0,+\infty)$ is a solution of \eqref{lineare1} if and only if
\[
z\in Z_0\cap Z_\infty.
\]
\end{remark}
\vs{12}
\ni
We conclude this subsection with some explicit formulas for solutions of the non-homogeneous equation
\begin{equation} \label{eqdicotomia}
Jz'+P(x)z=f,
\end{equation}
where $z, f\in L^2(0,1)$. They are based on the fact that the homogeneous equation \eqref{lineare1} has a suitable dichotomy at zero when $\l=0$ (see \cite{Co-book-78}).
\vs{12}
\ni
First of all, let us observe that from \cite[\textsection 3]{Co-book-78} we deduce that \eqref{eqdicotomia} has a solution $z_f\in L^\infty(0,1)$ when $f\in L^2(0,1)$. Moreover, let us point out that the previous results on the asymptotic behaviour for $x\to 0^+$ of the solutions of \eqref{lineare1} hold true also when $\l=0$; indeed, they are based on the fact that $\Delta^*>0$.
Hence, according to Remark \ref{spazio-soluzioni}, all the solutions $z\in L^2(0,1)$ of \eqref{eqdicotomia} are of the form
\[
z=c\zeta_{1,0}+z_f,
\]
for some $c\in \R$. More precisely, we have the following result:
\vs{12}
\ni
\begin{Theorem}\label{dicotomia-zero} (\cite[\textsection 3]{Co-book-78}) Let us consider $f\in L^2(0,1)$ and let $z\in L^2(0,1)$ be a solution of \eqref{eqdicotomia}. Then, there exist $c\in \R$ and $G:(0,1)\times (0,1)\lra \RR$ such that
\begin{equation} \label{formula-zero}
z(x)=c\zeta_{1,0}(x)+\int_0^1 G(x,\xi)f(\xi)\,d\xi,\quad \forall \ x\in (0,1).
\end{equation}
Moreover, there exist $K>0$ such that
\begin{equation} \label{G-limitato}
||G(x,\xi)||\leq
\begin{cases}
K \left(\dfrac{\min (x,\xi)}{\max (x,\xi)}\right)^{\sqrt{\Delta^{*}}}                                                                       & \text{if } \beta = 1 \vspace{4pt} \\
K \left( \dfrac{ e^{ -\min( x, \xi )^{ 1 - \beta } } }{ e^{ -\max( x, \xi )^{ 1 - \beta } } } \right)^{ \sqrt{ \Delta^{*} }/ ( \beta - 1 ) } & \text{if } \beta > 1
\end{cases}
\quad \leq K,\quad \forall \ (x,\xi)\in (0,1)\times (0,1).
\end{equation}
\end{Theorem}
\vs{8}
\ni
{\sc Proof.}
We just point out that the result follows from the change of variables $ x = \phi_{\beta}(t) $, from estimates in \cite[\textsection 3, formulas (3) and (4)]{Co-book-78} and Propositions \ref{stimezero2-0} or \ref{stimezero2}.
\qed
\subsection{Oscillatory properties} \label{oscillazione}
\vs{12}
\ni
In this subsection we develop an oscillatory theory for nontrivial solutions of \eqref{lineare1}, based on the study of the angular coordinate in the phase-plane (see \cite{We-book-87}). For every nontrivial solution $(u,v,\l)$ of \eqref{lineare1} let us introduce the polar coordinates $(\rho,\theta)=(\rho(x,\l),\theta(x,\l))$ according to
\[
\left\{\begin{array}{l}
			u=\rho \cos \theta\\
			\\
			v=\rho \sin \theta.
			\end{array}
\right.
\]
Observe that $\theta$ is defined mod. $2\pi$; we do not impose a normalization condition on $\theta$ and then the following results hold true for any angular coordinate associated to a nontrivial solution $z$. As a first step, we are able to study the asymptotic behaviour of $\theta$ when $x\to +\infty$ or $x\to 0^+$; this follows from the results of Subsection \ref{stimeasintotiche}.
\vs{12}
\ni
\begin{proposition} \label{limitethetainf} (\cite[Prop. 2.4]{CaDa-10})
For every $\l \in \Lambda$ the function $\theta(\cdot,\l)$ has limit at infinity and we have either
\begin{equation} \label{limitetheta1}
\lim_{x\to +\infty} \theta(x,\l)=\pi-\arctan \sqrt{\dfrac{\l-\mu^-}{\mu^+-\l}} \quad ({\rm{mod}} \ \pi)
\end{equation}
or
\begin{equation} \label{limitetheta2}
\lim_{x\to +\infty} \theta(x,\l)=\arctan \sqrt{\dfrac{\l-\mu^-}{\mu^+-\l}} \quad ({\rm{mod}} \ \pi).
\end{equation}
Moreover, \eqref{limitetheta1} and \eqref{limitetheta2} correspond to the cases when \eqref{limiteuv1} and \eqref{limiteuv2} are fulfilled, respectively.
\end{proposition}
\vs{12}
\ni
\begin{proposition} \label{limitethetazero}
For every $\l \in \Lambda$ the function $\theta(\cdot,\l)$ has limit at zero and we have either
\begin{equation} \label{limitethetazero1}
\lim_{x\to 0^+} \theta(x,\l)=\arctan\dfrac{w^*_{1,2}}{w^*_{1,1}} \quad ({\rm{mod}} \ \pi)
\end{equation}
or
\begin{equation} \label{limitethetazero2}
\lim_{x\to 0^+} \theta(x,\l)=\arctan \dfrac{w^*_{2,2}}{w^*_{2,1}} \quad ({\rm{mod}} \ \pi),
\end{equation}
where $w^*_1$ and $w^*_2$ are eigenvectors of $C$ associated to the eigenvalues $\sigma^{\pm}$, respectively. Moreover, \eqref{limitethetazero1} and \eqref{limitethetazero2} correspond to the cases when \eqref{limiteuvzero1} and \eqref{limiteuvzero2} are fulfilled, respectively.
\end{proposition}
\vs{12}
\ni
Let us observe that the possible limits of $\theta (\cdot,\l)$ at zero do not depend on $\l\in\Lambda$; in what follows, we denote
\[
\begin{array}{l}
\theta (+\infty,\l)=\lim_{x\to +\infty} \theta(x,\l)\\
\\
\theta (0)=\lim_{x\to 0^+} \theta(x,\l),
\end{array}
\]
which exist and are finite by Proposition \ref{limitethetainf} and Proposition \ref{limitethetazero}.
\vs{12}
\ni
\begin{remark} \label{theta-soluzioni} According to Remark \ref{spazio-soluzioni} and the above Propositions, we deduce that
for a nontrivial solution $z$ of \eqref{lineare1} we have
\[
z\in L^2(1,+\infty) \quad \Longleftrightarrow \quad \theta (+\infty,\l)=\pi-\arctan \sqrt{\dfrac{\l-\mu^-}{\mu^+-\l}} \quad ({\rm{mod}} \ \pi)
\]
and
\[
z\in L^2(0,1) \quad \Longleftrightarrow \quad \theta (0)=\arctan\dfrac{w^*_{1,2}}{w^*_{1,1}} \quad ({\rm{mod}} \ \pi).
\]
\end{remark}
\vs{12}
\ni
Proposition \ref{limitethetainf} and Proposition \ref{limitethetazero} imply that any angular function $\theta (\cdot,\l)$ is bounded on $(0,+\infty)$, for every $\l\in \Lambda$. As a consequence, we can associate to every nontrivial solution $z$ of \eqref{lineare1} the rotation number
\begin{equation} \label{defrotazioni}
\rot \ (z)=\dfrac{\theta(+\infty,\l)-\theta(0)}{\pi}.
\end{equation}
Roughly speaking, the unboundedness of the interval and the singularity at zero do not prevent solutions to perform only a finite number of rotations around the origin (as in the regular case). It is important to observe that $\rot \ (z)$ does not depend on the choice of the angular function of $z$. In Section \ref{sez-nonlineare} we will study some continuity properties of the rotation number defined in \eqref{defrotazioni}.
\vs{12}
\ni
We conclude this subsection with some asymptotic phase-plane analysis for \eqref{lineare1}; as above, we prove the results for $x\to +\infty$. The case of $x\to 0^+$ can be obtained in an analogous way by means of the change of variable $ x = \phi_{ \beta }( t ) $ already introduced.
\vs{6}
\ni
Let us consider again \eqref{sistemaintermedio}, which is equivalent to \eqref{lineare1}, and a similar system
\begin{equation} \label{lineare2}
Jz'+\tilde P(x)z=\tilde \l z,
\end{equation}
where $\tilde P\in {\cal P}_{\mu} $ and $\tilde \l\in \Lambda$; \eqref{lineare2} can be written in the form
\begin{equation} \label{intermedio2}
z'=B_{\tilde \l}z+\tilde Q(x)z,
\end{equation}
where $B_{\tilde \l}=J^{-1}(\tilde \l \mbox{Id }-P_\infty)$ and $\tilde Q(x)=J^{-1}(P_\infty-\tilde P(x))$, for every $x>0$. Let us note that the matrix $P_\infty$ is the same both for $P$ and $\tilde P$, since $P, \tilde P\in  {\cal P}_{\mu}$.
\vs{4}
\ni
For every $\l\in \Lambda$, let $b_{1,\l}$, $b_{2,\l}$ be as in Proposition \ref{stimeinf2}; from the discussion leading to Proposition \ref{stimeinf2} we know that
\[
b_{1,\l}=(\l-\mu^+,\sqrt{\Delta_\l}),\quad b_{2,\l}=(-\l+\mu^+,\sqrt{\Delta_\l}),
\]
for every $\l \in \Lambda$; moreover, there exists $\rho_\l>0$ such that
\[
\begin{array}{l}
b_{1,\l}=\rho_\l (\cos \theta_{\infty,\l},\sin \theta_{\infty,\l})\\
\\
b_{2,\l}=\rho_\l (-\cos \theta_{\infty,\l},\sin \theta_{\infty,\l}),
\end{array}
\]
where
\[
\theta_{\infty,\l}=\pi-\arctan \sqrt{\dfrac{\l-\mu^-}{\mu^+-\l}} \quad ({\rm{mod}} \ \pi).
\]
For every $\theta \in \R$, let $r_\theta$ be the straight line of equation $x\sin \theta-y\cos \theta=0$ and let $v_\theta=(\sin \theta,-\cos \theta)$; moreover, let $r^{\pm}_\theta$ be the half-lines given by the intersection of $r_\theta$ with the half-planes $H^+=\{(x,y)\in \RR:\ x>0\}$ and $H^-=\{(x,y)\in \RR:\ x<0\}$, respectively. We are in position to prove the following result:
\vs{12}
\ni
\begin{proposition} \label{campo-vettoriale-1} For every $\tilde \l\in \Lambda$, $\tilde P\in {\cal P}_{\mu}$ and for every $\theta \in (\pi/2,\pi)$ there exist $\tilde \delta>0$ and  $\tilde x_{\infty}=\tilde x_{\infty}(\tilde \l,\tilde P,\theta)>0$ such that for every $\l\in \Lambda$ and $P\in  {\cal P}_{\mu}$ with
\begin{equation} \label{coefficienti-vicini}
|\l-\tilde \l|<\tilde \delta,\quad ||P-\tilde P||_{L^\infty(1,+\infty)}<\tilde \delta
\end{equation}
we have
\begin{equation} \label{direzione-campo-inf-1}
\begin{array}{l}
\theta<\theta_{\infty,\l}\quad \Rightarrow \quad \langle v_\theta, B_\l w+Q(x)w \rangle > 0,\quad \forall \ w\in r^-_\theta,\ \forall\ x\geq \tilde x_{\infty}\\
\\
\theta>\theta_{\infty,\l}\quad \Rightarrow \quad \langle v_\theta, B_\l w+Q(x)w \rangle < 0,\quad \forall \ w\in r^-_\theta,\ \forall\ x\geq \tilde x_{\infty}.
\end{array}
\end{equation}
\end{proposition}
\vs{8}
\ni
{\sc Proof.} First of all, let us observe that it is sufficient to prove \eqref{direzione-campo-inf-1} when $w$ is a versor. Therefore, let  $w=(\cos \theta,\sin \theta)$; a simple computation shows that
\begin{equation} \label{conto1}
\phi_\l(\theta,w):= \langle v_\theta, B_\l w \rangle =\cos^2 \theta \left((\mu^+-\l)\tan^2 \theta -(\l-\mu^-)\right),\quad \forall \ \l\in \Lambda.
\end{equation}
Let us fix $\tilde \l\in \Lambda$, $\tilde P\in {\cal P}_{\mu}$ and $\theta \in (\pi/2,\pi)$ such that $\theta<\theta_{\infty,\tilde \l}$; the continuity of $\theta_{\infty,\l}$ as a function of $\l\in \Lambda$ implies that there exists $\delta_1>0$ such that $\theta<\theta_{\infty,\l}$ if $|\l-\tilde \l|<\delta_1$.
\vs{4}
\ni
From \eqref{conto1} we deduce that
\[
\begin{array}{l}
\phi_{\tilde \l}(\theta_{\infty,\tilde \l},w)=0\\
\\
\theta<\theta_{\infty,\tilde \l}\quad \Rightarrow \quad \phi_{\tilde \l}(\theta,w) >\phi_{\tilde \l}(\theta_{\infty,\tilde \l},w)=0.
\end{array}
\]
Hence, there exists $\delta_2\in (0,\delta_1)$ such that
\begin{equation} \label{stima-oggi-1}
|\l-\tilde \l|<\delta_2 \quad \Rightarrow \quad
\theta<\theta_{\infty,\l}\quad {\mbox{and}} \quad \phi_\l(\theta,w) >\dfrac{\phi_{\tilde \l}(\theta,w)}{4}>0
\end{equation}
Now, from assumption \eqref{limitiPinfinito} we deduce that
\[
\lim_{x\to +\infty} \langle v_\theta, \tilde Q(x)w \rangle =\lim_{x\to +\infty} \langle v_\theta, J^{-1}(P_\infty-\tilde P(x))w \rangle =0;
\]
this implies that there exists $\tilde x_\infty=\tilde x_\infty (\tilde \l,\tilde P,\theta)>1$ such that
\begin{equation} \label{stima-oggi-2}
x\geq \tilde x_\infty \quad \Rightarrow \quad | \langle v_\theta, \tilde Q(x)w \rangle | < \dfrac{\phi_{\tilde \l}(\theta,w)}{16}.
\end{equation}
On the other hand, setting $\delta_3=\phi_{\tilde \l}(\theta,w)/16$, if $||P-\tilde P||_{L^\infty(1,+\infty)}<\delta_3$ we have
\begin{equation} \label{stima-oggi-3}
\begin{array}{ll}
| \langle v_\theta, Q(x)w \rangle - \langle v_\theta, \tilde Q(x)w \rangle | = | \langle v_\theta,J^{-1}(\tilde P(x)-P(x))w \rangle |\leq \\
\\
\leq ||\tilde P(x)-P(x)||< \dfrac{\phi_{\tilde \l}(\theta,w)}{16},\quad \forall \ x\geq 1.
\end{array}
\end{equation}
From \eqref{stima-oggi-2} and \eqref{stima-oggi-3} we deduce that
\begin{equation} \label{stima-oggi-4}
||P-\tilde P||_{L^\infty(1,+\infty)}< \delta_3,\quad x\geq \tilde x_\infty \quad \Rightarrow \quad | \langle v_\theta, Q(x)w \rangle |<\dfrac{\phi_{\tilde \l}(\theta,w)}{8}.
\end{equation}
Now, let us set $\tilde \delta=\min (\delta_2,\delta_3)$; when $|\l-\tilde \l|<\tilde \delta$ and $||P-\tilde P||_{L^\infty(1,+\infty)}<\tilde \delta$ both \eqref{stima-oggi-1} and \eqref{stima-oggi-4} hold true. As a consequence, we obtain
\begin{equation} \label{stima-oggi-5}
x\geq \tilde x_\infty \quad \Rightarrow \quad \langle v_\theta, B_\l w+Q(x)w \rangle > \dfrac{\phi_{\tilde \l}(\theta,w)}{4}-\dfrac{\phi_{\tilde \l}(\theta,w)}{8}>0,
\end{equation}
i.e. the first inequality in \eqref{direzione-campo-inf-1} is satisfied.

\ni
An analogous argument proves the validity of the second inequality in \eqref{direzione-campo-inf-1}.
\qed
\vs{12}
\ni
In a very similar way it is possible to prove the following Proposition:
\vs{12}
\ni
\begin{proposition} \label{campo-vettoriale-2}
For every $\tilde \l\in \Lambda$, $\tilde P\in {\cal P}_{\mu} $ and for every $\theta \in (0,\pi/2)$ there exist $\tilde \delta_1>0$ and  $\tilde x_{\infty,1}=\tilde x_{\infty,1}(\tilde \l,\tilde P,\theta)>0$ such that for every $\l\in \Lambda$ and $P\in {\cal P}_{\mu} $ with
\begin{equation} \label{coefficienti-vicini-2}
|\l-\tilde \l|<\tilde \delta,\quad ||P-\tilde P||_{L^\infty(1,+\infty)}<\tilde  \delta
\end{equation}
we have
\begin{equation} \label{direzione-campo-inf-2}
\begin{array}{l}
\theta<\pi-\theta_{\infty,\l}\quad \Rightarrow \quad \langle v_\theta, B_\l w+Q(x)w \rangle <0,\quad \forall \ w\in r^+_\theta,\ \forall\ x\geq \tilde x_{\infty,1}\\
\\
\theta>\pi-\theta_{\infty,\l}\quad \Rightarrow \quad \langle v_\theta, B_\l w+Q(x)w \rangle >0,\quad \forall \ w\in r^+_\theta,\ \forall\ x\geq \tilde x_{\infty,1}.
\end{array}
\end{equation}
\end{proposition}
\vs{12}
\ni
From Proposition \ref{campo-vettoriale-1} and Proposition \ref{campo-vettoriale-2} we deduce the following result:
\vs{12}
\ni
\begin{proposition} \label{confini-theta-inf} For every $\tilde \l\in \Lambda$ and $\tilde P\in {\cal P}_{\mu}$ there exists $\tilde \epsilon>0$ such that for every $\epsilon\in (0,\tilde \epsilon)$ there exist $\tilde \delta>0$ and $\tilde x_{\infty}=\tilde x_{\infty}(\tilde \l,\tilde P,\epsilon)>0$ such that for every $\l\in \Lambda$ and $P\in {\cal P}_{\mu}$ with
\begin{equation} \label{coefficienti-vicini-3}
|\l-\tilde \l|<\tilde \delta,\quad ||P-\tilde P||_{L^\infty(1,+\infty)}<\tilde \delta
\end{equation}
and for every nontrivial solution $z \in L^{2}( (1, +\infty ) ) $ of \eqref{lineare1} we have
\begin{equation} \label{confine-inf}
|\theta (x,\l)-\theta(+\infty,\l)|<\epsilon, \quad \forall \ x\geq \tilde x_{\infty},
\end{equation}
where $\theta (\cdot,\l)$ is any angular coordinate of $z$.
\end{proposition}
\vs{8}
\ni
{\sc Proof.} 
Without loss of generality let us assume that
\[
\theta_{ \infty, \lambda } = \pi - \arctan \sqrt{ \dfrac{ \lambda - \mu^{-} }{ \mu^{+} - \lambda } } \in \left( \frac{ \pi }{ 2 }, \pi \right)
\]
and define
\[
\tilde \epsilon = \min\left\{ \theta_{ \infty, \tilde \lambda } - \frac{ \pi }{ 2 },
                              \pi - \theta_{ \infty, \tilde \lambda } \right\}
                > 0.
\]
Fix any $ \epsilon \in ( 0, \tilde \epsilon ) $ and consider
\[
\theta_{1} = \theta_{ \infty, \tilde \lambda } - \frac{ \epsilon }{ 2 } \qquad \text{and} \qquad
\theta_{2} = \theta_{ \infty, \tilde \lambda } + \frac{ \epsilon }{ 2 },
\]
thus the cone between  $ r_{ \theta_{1} } $ and $ r_{ \theta_{2} } $ lies inside the II and the IV quadrants and its angular amplitude is exactly $ \epsilon $.
We use the continuity of $ \theta_{ \infty, \lambda } $ with respect to $ \lambda $ and apply Proposition~\ref{campo-vettoriale-1} twice
with the choices $ \theta = \theta_{1} $ and $ \theta = \theta_{2} $ in order to find $ \tilde \delta > 0 $ and
$ \tilde x_{\infty} = \tilde x_{\infty}( \tilde \lambda, \tilde P, \epsilon ) > 0 $ in such a way that, if
\eqref{coefficienti-vicini-3} hold, then
$ \left| \theta_{ \infty, \lambda } - \theta_{ \infty, \tilde \lambda } \right| < \epsilon / 2 $ and \eqref{direzione-campo-inf-1} hold.
We remark that $ \tilde x_{\infty} $ depends only on $ \tilde \lambda, \tilde P, \epsilon $ since the number $ \tilde x_{\infty,1} $ provided by Proposition~\ref{campo-vettoriale-1} depends on $ \theta_{1} $ and $ \theta_{2} $ which depend only on $ \tilde \lambda $ and $ \epsilon $.

By construction we have $ \theta_{1} < \theta_{ \infty, \lambda } < \theta_{2} $ and \eqref{direzione-campo-inf-1} implies that the vector field of
\eqref{sistemaintermedio} points strictly outwards the cone between $ r_{ \theta_{1} } $ and $ r_{ \theta_{2} } $ for all $ x \ge \tilde x_{\infty} $.
Therefore, any nontrivial solution $ z \in L^{2}( 1, +\infty ) $ of \eqref{sistemaintermedio} 
approaches the origin at the angle $ \theta_{ \infty, \lambda } $  as $ x $ tends to infinity and
a standard phase plane argument shows that $ z(x) $ must remain inside the cone between
$ r_{ \theta_{1} } $ and $ r_{ \theta_{2} } $ for all $ x \ge \tilde x_{\infty} $.
Hence \eqref{confine-inf} follows.
\qed

\vs{12}
\ni
By means of the transformation $x=\phi_\beta (t)$, it is possible to prove an analogous result concerning the local behaviour of the angular coordinate when $x\to 0^+$; indeed, we have the following:
\vs{12}
\ni
\begin{proposition} \label{confini-theta-zero} There exist $\epsilon_0>0$ such that for every $\epsilon\in (0,\epsilon_0)$, $\tilde \l\in \Lambda$ and $\tilde P\in  {\cal P}_{\mu}$ there exist $\tilde \delta_0>0$ and $x_0=x_0(\tilde P, \epsilon)>0$ such that for every $\l\in \Lambda$ and $P\in {\cal P}_{\mu}$ with
\begin{equation} \label{coefficienti-vicini-4}
|\l-\tilde \l|<\tilde \delta_0,\quad ||P-\tilde P||_{L^\infty(0,1)}<\tilde \delta_0
\end{equation}
and for every nontrivial solution $z \in L^{2}( (0, 1 ) )$ of \eqref{lineare1} we have
\begin{equation} \label{confine-zero}
|\theta (x,\l)-\theta(0)|<\epsilon, \quad \forall \ x\in (0,x_0],
\end{equation}
where $\theta (\cdot,\l)$ is any angular coordinate of $z$.
\end{proposition}
\vs{12}
\ni
\section{The linear eigenvalue problem} \label{sez-autovalori}
\vs{12}
\ni
In this Section we are dealing with the study of the spectral theory for the linear operator formally defined by
\begin{equation} \label{tau}
\tau z=Jz'+P(x)z,\ x>0,
\end{equation}
where $P\in {\cal P}_{\mu}$. Some information on the spectrum of $\tau$ follow directly from a standard spectral theory (see e.g. \cite{ScTr-02,We-book-87}). Indeed, \cite[Th. 6.8]{We-book-87} ensures that $\tau$ is in the limit point case at infinity; moreover, from Remark \ref{spazio-soluzioni} we deduce that $\tau$ is in the limit point case also at zero. Let us point out that this fact is a consequence of assumption $({\cal P}_3)$ on $P^*$.
\vs{6}
\ni
Let us consider the operator $A_0$ defined by
\begin{equation} \label{def-A0}
\begin{array}{l}
D(A_0)=\{z\in L^2(0,+\infty):\ z\in AC(0,+\infty),\ \tau z\in L^2(0,+\infty)\},\\
\\
A_0 z=\tau z,\quad \forall \ z\in D(A_0).
\end{array}
\end{equation}
From \cite[Th. 5.8]{We-book-87} we deduce that $A_0$ is the unique self-adjoint realization of $\tau$; moreover, arguing as in the proof of \cite[Lemma 5.1]{ScTr-02}, it is possible to see that $\sigma_{\rm{ess}} (A_0)=(-\infty,\mu^-]\cup [\mu^+,+\infty)$ .
\vs{4}
\ni
As far as $D_0:=D(A_0)$ is concerned, we are able to prove the following result:
\vs{12}
\ni
\begin{proposition} \label{dominio}
For every $z\in D_0$ we have
\[
z\in H^1 (1,+\infty),\quad z\in L^\infty (0,+\infty).
\]
\end{proposition}
\vs{8}
\ni
{\sc Proof.}
Assume that $z\in D_0$. Since $P\in L^\infty(1,+\infty)$ we deduce that $P(x)z\in L^2(1,+\infty)$; hence $Jz'=\tau z-P(x)z\in L^2 (1,+\infty)$. This proves that $z\in H^1 (1,+\infty)\subset L^\infty (1,+\infty)$.
\vs{4}
\ni
The fact that $z\in L^\infty (0,1)$ immediately follows from \eqref{formula-zero} and \eqref{G-limitato}.
\qed
\vs{12}
\ni
The aim of this Section is to study the problem of the existence of eigenvalues of $A_0$ in $\Lambda$; first of all, let us observe that every eigenvalue of $A_0$ is simple, since $\tau$ is in the limit point case at infinity. Moreover, from Remark \ref{spazio-soluzioni} we know that $\l\in \Lambda$ is an eigenvalue of $A_0$ if and only if there exists $c_\l\in \R$ such that
\begin{equation} \label{car-aut-sol-fond}
\zeta_{1,\l}=c_\l z_{1,\l},
\end{equation}
where $z_{1,\l}$ and $\zeta_{1,\l}$ are given in Proposition \ref{stimeinf2} and Proposition \ref{stimezero2}, respectively.
\vs{12}
\ni
\begin{remark} \label{regolarita-sol} According to Lemma \ref{limitiuv} and Lemma \ref{limitiuvzero}, when $\l\in \Lambda$ is an eigenvalue of $A_0$ the associated eigenfunction $z_\l$  satisfies $z_{\l}\in H^1_0(0,+\infty)$.
\end{remark}
\vs{12}
\ni
In what follows we show that it is possible to write a condition equivalent to \eqref{car-aut-sol-fond} by means of the angular function $\theta$ associated to solutions of \eqref{lineare1} introduced in Subsection \ref{oscillazione}. To this aim, let us denote by $\vartheta (\cdot,\l)$ the angular coordinate of $\zeta_{1,\l}$, normalized in such a way that $\vartheta (0)\in (0,\pi)$, for every $\l\in \Lambda$.
\vs{4}
\ni
From Proposition \ref{limitethetainf} we know that there exists
\[
\lim_{x\to +\infty} \vartheta (x,\l)=\vartheta (+\infty,\l)
\]
and that this limit corresponds to a function belonging to $H^1(1,+\infty)$ if and only if
\begin{equation} \label{limite-autovalori}
\vartheta (+\infty,\l)=\pi -\arctan \sqrt{\dfrac{\l-\mu^-}{\mu^+-\l}} \quad ({\rm{mod}} \ \pi).
\end{equation}
Let us define $\nu:\Lambda \to \R$ by
\[
\nu (\l)=\lim_{x\to +\infty} \vartheta (x,\l),\quad \forall \ \l\in \Lambda.
\]
We then have the following characterization of the eigenvalues of $A_0$:
\vs{12}
\ni
\begin{Theorem} \label{condizione-autovalore} A number $\l \in \Lambda$ is an eigenvalue of $A_0$ if and only if
\begin{equation} \label{limite-autovalore-2}
\nu (\l)=\pi -\arctan \sqrt{\dfrac{\l-\mu^-}{\mu^+-\l}} \quad ({\rm{mod}} \ \pi).
\end{equation}
\end{Theorem}
\vs{12}
\ni
In order to prove the existence of eigenvalues of $A_0$ it is then sufficient to study the behaviour of the function $\nu^*:\Lambda \to {\R}$ defined by
\[
\nu^*(\l)=\nu(\l)+\arctan \sqrt{\dfrac{\l-\mu^-}{\mu^+-\l}},\quad \forall \ \l \in \Lambda.
\]
We will prove that $\nu^*$ is strictly increasing and continuous in $\Lambda$.
\vs{12}
\ni
\begin{proposition} \label{monotonia-nu*} The function $\nu^*:\Lambda \to {\R}$ is strictly increasing in $\Lambda$.
\end{proposition}
\vs{8}
\ni
{\sc Proof.} Let us first observe that $\nu^*$ is the sum of $\nu$ and of the function $\nu_*$ defined by
\[
\nu_*(\l)=\arctan \sqrt{\dfrac{\l-\mu^-}{\mu^+-\l}},\quad \forall \ \l\in \Lambda;
\]
since $\nu_*$ is strictly increasing in $\Lambda$, it is sufficient to prove that $\nu$ is increasing in $\Lambda$.
\vs{4}
\ni
To this aim, let us recall (cf. \cite[Cor. 16.2]{We-book-87}) that for every fixed $x>0$ the function
\[
\begin{array}{ll}
\varphi_x:&\Lambda \to \R\\
&\\
&\l\to \vartheta (x,\l)
\end{array}
\]
is increasing in $\Lambda$.
\vs{4}
\ni
Now, let $\l, \l'\in \Lambda$ with $\l<\l'$; for every $x>0$ we have
\[
\vartheta (x,\l)\leq \vartheta (x,\l');
\]
passing to the limit for $x\to +\infty$ we obtain
\[
\lim_{x\to +\infty} \vartheta (x,\l)\leq \lim_{x\to +\infty} \vartheta (x,\l'),
\]
i.e.
\[
\nu(\l)\leq \nu (\l').
\]
\qed
\vs{12}
\ni
\begin{proposition} \label{continuita-nu*} The function $\nu^*:\Lambda \to {\R}$ is continuous.
\end{proposition}
\vs{8}
\ni
{\sc Proof.} Let us observe again that it is sufficient to prove the continuity of $\nu$. To this aim, let us fix $\tilde \l\in \Lambda$;
let us consider $\epsilon>0$ sufficiently small and apply Proposition \ref{confini-theta-inf} and Proposition \ref{confini-theta-zero} with $\tilde P=P$. Let $\delta_1=\min (\tilde \delta,\tilde \delta_0)$ and let us denote by $x_\infty$ and $x_0$ the numbers given in those Propositions.
\vs{4}
\ni
Let us recall that a usual continuous dependence argument on the interval $[x_0,x_\infty]$, on which the equation \eqref{lineare1} is not singular, proves that there exists $\delta_2>0$ such that if $|\l-\tilde \l|<\delta_2$ then
\begin{equation} \label{rotazioni-non-singolare}
|(\vartheta (x_\infty, \l)-\vartheta (x_0, \l))-(\vartheta (x_\infty,\tilde \l)-\vartheta (x_0,\tilde \l))|<\epsilon.
\end{equation}
Consider now $\delta=\min (\delta_1,\delta_2)$ and assume that $|\l-\tilde \l|<\delta$; we can write
\begin{equation} \label{nuovo-1}
\begin{array}{l}
\nu(\l)-\nu(\tilde \l)=\vartheta (+\infty,\l)-\vartheta (+\infty,\tilde \l)=\vartheta (+\infty,\l)-\vartheta (x_\infty,\l)+\\
\\
+\vartheta (x_\infty,\l)-\vartheta (x_0,\l)+\vartheta (x_0,\tilde \l)-\vartheta (x_\infty, \tilde \l)+\vartheta (x_0,\l)-\vartheta (0,\l)+\\
\\
-\vartheta (x_0,\tilde \l)+\vartheta (0,\tilde \l)+\vartheta(x_\infty,\tilde \l)-\vartheta(+\infty,\tilde \l),
\end{array}
\end{equation}
taking into account that $\vartheta (0,\l)=\vartheta (0,\tilde \l)$. From Proposition \ref{confini-theta-inf} and Proposition \ref{confini-theta-zero} we deduce that
\begin{equation} \label{nuovo-2}
\begin{array}{l}
|\vartheta (+\infty,\l)-\vartheta (x_\infty,\l)|<\epsilon,\quad |\vartheta(x_\infty,\tilde \l)-\vartheta(+\infty,\tilde \l)|<\epsilon\\
\\
|\vartheta (x_0,\l)-\vartheta (0,\l)|<\epsilon,\quad |\vartheta(x_0,\tilde \l)-\vartheta(0,\tilde \l)|<\epsilon.
\end{array}
\end{equation}
From \eqref{rotazioni-non-singolare}, \eqref{nuovo-1} and \eqref{nuovo-2} we obtain
\[
|\nu(\l)-\nu(\tilde \l)|<5\epsilon
\]
and this concludes the proof.
\qed
\vs{12}
\ni
For every $k\in {\Z}$, let us denote by $\l_k\in \Lambda$ (if it exists) the number such that
\[
\nu^*(\l_k)=k\pi,
\]
i.e.
\begin{equation} \label{autovalore-k}
\vartheta (+\infty,\l_k)=k\pi+\pi-\arctan \sqrt{\dfrac{\l_k-\mu^-}{\mu^+-\l_k}}.
\end{equation}
The number $\l_k$ is the '$k$-th eigenvalue' of $A_0$ (if it exists) and we denote by $z_k\in D_0$ the corresponding eigenfunction; recalling \eqref{limitethetazero1}, \eqref{defrotazioni}, from \eqref{autovalore-k} and the fact that $\nu^*$ is strictly increasing we immediately deduce the following result:
\vs{12}
\ni
\begin{proposition} \label{numero-rotazioni-autofunzioni} For every $k\in {\Z}$ we have
\begin{equation} \label{num-rot-aut}
\begin{array}{ll}
{\mbox{rot $(z_k)$}}\in (k,k+1)&\displaystyle{{\mbox{if $\arctan \dfrac{w^*_{1,2}}{w^*_{1,1}} \in (0,\pi/2)$}}}\\
&\\
{\mbox{rot $(z_k)$}}\in (k-1/2,k+1/2)&\displaystyle{{\mbox{if $\arctan \dfrac{w^*_{1,2}}{w^*_{1,1}} \in (\pi/2,\pi)$}}.}
\end{array}
\end{equation}
Moreover, for every $k, l\in {\Z}$ with $k\neq l$ we also have
\begin{equation} \label{indici-diversi}
{\mbox{rot $(z_k)$}}\neq {\mbox{rot $(z_l)$}}.
\end{equation}
\end{proposition}
\vs{12}
\ni
In what follows, we give some results on the accumulation of eigenvalues of $A_0$ at the boundary of $\Lambda$. We consider the (possible) accumulation at the end-point $\mu^+$; conditions for accumulation at $\mu^-$ can be obtained in an analogous way.
\vs{4}
\ni
From \eqref{limite-autovalore-2} and the definition of $\nu^*$ we infer that the existence of eigenvalues accumulating at $\mu^+$ depends on the behaviour of $\nu^*$ in a left neighbourhood of $\mu^+$. This behaviour can be described by means of the limit
\begin{equation} \label{limite-nu*}
\lim_{\l \to (\mu^+)^-} \nu^* (\l),
\end{equation}
whose existence is guaranteed from Proposition \ref{monotonia-nu*}; more precisely, when the limit in \eqref{limite-nu*} is infinite, then there exists $k_0\in \Z$ such that for every $k\in \Z$, $k\geq k_0$, there exists $\l_k \in \Lambda$ for which \eqref{autovalore-k} holds true and
\[
\lim_{k\to +\infty} \l_k=\mu^+,
\]
i.e. there is accumulation of eigenvalues at $\mu^+$. On the other hand, when the limit in \eqref{limite-nu*} is finite, then there exists $M^+\in \R$ such that
\[
\nu^*\left(\dfrac{\mu^++\mu^-}{2}\right)<\nu^*(\l)<M^+,\quad \forall \ \l\in \left(\dfrac{\mu^++\mu^-}{2},\mu^+\right);
\]
this implies that there is at most a finite number of eigenvalues of $A_0$ in $((\mu^++\mu^-)/2,\mu^+)$, i.e. there is not accumulation of eigenvalues at $\mu^+$.
\vs{6}
\ni
Now, let us observe that the fact that the limit in \eqref{limite-nu*} is finite or infinite depends on the analogous limit
\begin{equation} \label{limite-nu}
\lim_{\l \to (\mu^+)^-} \nu (\l),
\end{equation}
since the function $\nu_*$ is bounded in $\Lambda$. We are able to show that the finiteness of the limit in \eqref{limite-nu} depends on the behaviour of \eqref{lineare1} when $\l=\mu^+$; to this aim, let us observe that a more careful analysis proves that Proposition \ref{stimezero2} and Proposition \ref{limitethetazero} hold true also when $\l=\mu^+$. This implies that we are allowed to consider the solution $\zeta_{1,\l}$ of \eqref{lineare1} with $\l=\mu^+$ satisfying \eqref{limitezeta1} and the corresponding angular coordinate $\vartheta (\cdot,\mu^+)$, normalized in such a way that $\vartheta (0,\mu^+)\in (0,\pi)$.
\vs{12}
\ni
\begin{lemma} \label{nu-finito} Assume
\begin{equation} \label{confrfinito-2}
\lim_{x\to +\infty} \vartheta (x,\mu^+)=\theta^+\in \R;
\end{equation}
then we have
\begin{equation} \label{no-acc-1}
\lim_{\l \to (\mu^+)^-} \nu(\l)<+\infty.
\end{equation}
\end{lemma}
\vs{8}
\ni
{\sc Proof.} Let us observe that \eqref{confrfinito-2} implies that there exist $\Phi\in \R$ and $X>0$ such that
\[
\vartheta (x,\mu^+)<\Phi,\quad \forall \ x\geq X.
\]
Moreover, from the monotonicity of $\vartheta (x,\cdot)$, for every $x\geq X$, we deduce that
\[
\vartheta (x,\l)\leq \vartheta (x,\mu^+)<\Phi,\quad \forall \  \l<\mu^+.
\]
Therefore, for every $\l\in \Lambda$ the function $\vartheta (\cdot ,\l)$ is bounded from above by $\Phi$ in $[X,+\infty)$, hence we have
\[
\nu(\l)=\lim_{x\to +\infty} \vartheta (x,\l)\leq \Phi,\quad \forall \ \l <\mu^+.
\]
This is sufficient to conclude that \eqref{no-acc-1} holds true.
\qed
\vs{12}
\ni
\begin{lemma} \label{nu-infinito} Assume
\begin{equation} \label{confrinfinito-2}
\lim_{x\to +\infty} \vartheta (x,\mu^+)=+\infty
\end{equation}
and that there exists $X>0$ such that
\begin{equation} \label{monotonia-debole}
p_{11}(x)<\mu^-,\quad \forall \ x\geq X.
\end{equation}
Then we have
\begin{equation} \label{acc-2}
\lim_{\l \to (\mu^+)^-} \nu(\l)=+\infty.
\end{equation}
\end{lemma}
\vs{8}
\ni
{\sc Proof.} Let us first observe that for every $\l\in \R$ the angular function $\vartheta (\cdot, \l)$ satisfies the differential equation
\begin{equation} \label{eq-theta}
\theta'=(\l-p_{11}(x))\cos^2 \theta-2p_{12}(x)\cos \theta \sin \theta +(\l-p_{22}(x))\sin^2 \theta.
\end{equation}
From \eqref{eq-theta} and \eqref{monotonia-debole} we deduce that
\[
\forall\  x\geq X,\  \l>\mu^-:\quad \vartheta (x,\l)=0 \ ({\rm{mod}} \ \pi)\quad \Rightarrow \quad \vartheta' (x,\l)>0;
\]
hence, if there exist $k\in \Z$ and $x_k\geq X$ such that
\[
\vartheta (x_k,\l)>k\pi,
\]
for some $\l>\mu^-$, then we can conclude that
\[
\vartheta (x,\l)>k\pi,\quad \forall \ x\geq x_k.
\]
\vs{4}
\ni
Now, let us note that \eqref{confrinfinito-2} implies that for every $M>0$ there exists $x_M\geq X$ such that
\[
\vartheta (x,\mu^+)>M+2+\pi,\quad \forall\  x\geq x_M
\]
and let us fix $X^+\geq x_M$; the continuity of $\vartheta (X^+,\cdot)$ ensures that there exists $\l_M<\mu^+$ such that
\[
\vartheta (X^+,\l)>M+1+\pi,\quad \forall \ \l \in (\l_M,\mu^+).
\]
According to the above remark, this implies that
\[
\vartheta (x,\l)>M+1,\quad \forall  \  x\geq X^+,\ \l\in (\l_M,\mu^+)
\]
and then
\begin{equation} \label{dimostrazione-1}
\nu(\l)=\lim_{x\to +\infty} \vartheta (x,\l)> M,\quad \forall \ \l\in (\l_M,\mu^+).
\end{equation}
Therefore, for every $M>0$ there exists $\l_M<\mu^+$ such that \eqref{dimostrazione-1} holds, i.e.
\[
\lim_{\l \to (\mu^+)^-} \nu(\l)=+\infty.
\]
\qed
\vs{12}
\ni
The question of the existence of eigenvalues can be dealt, arguing as in the proof of Proposition 3.18 in \cite{CaDa-10}, as follows.
\vs{12}
\ni
\begin{proposition} \label{esistenza-autovalori} Assume that $P$ has the form \eqref{Dirac-magn}, where $\mu_a\in {\R}$, $k\in {\bf Z}\setminus \{0\}$ and $V\in C^1(0,+\infty)$ is a strictly increasing negative potential satisfying \eqref{ipotesi-potenziale-infinito}, with $\gamma_\infty<0$ and $\alpha_\infty \in (0,1]$, and \eqref{ipotesi-potenziale-zero}.

\ni
Then, the selfadjoint extension $A_0$ of the corresponding operator $\tau$ has a sequence of eigenvalues in $(-1,1)$ accumulating at $\lambda=1$.
\end{proposition}
\vs{8}
\ni
{\sc Proof.} We follow the same argument of \cite[Prop. 3.15]{CaDa-10}. We first observe that the differential equation satisfied by $\vartheta(\cdot,1)$ is
\[
\vartheta'(x,1)=1-\langle Q_{P(x)}[\cos \vartheta,\sin \vartheta],[\cos \vartheta,\sin \vartheta]\rangle,
\]
where $Q_{P(x)}$ denotes the quadratic form associated to the matrix $P(x)$. By computing the eigenvalues of $P(x)$, we can prove that
\[
\vartheta'(x,1)\geq 1-V(x)-\sqrt{1+\left(\dfrac{k}{x}+\mu_a V'(x)\right)^2},\quad \forall \ x\geq 1.
\]
From assumption \eqref{ipotesi-potenziale-infinito} we infer that 
\[
1-V(x)-\sqrt{1+\left(\dfrac{k}{x}+\mu_a V'(x)\right)}=-\dfrac{\gamma_\infty}{x^\alpha}+o\left(\dfrac{1}{x^\alpha}\right),\quad x\to +\infty;
\]
this is sufficient to conclude that
\[
\lim_{x\to +\infty} \vartheta(x,1)=+\infty.
\]
The result then follows from the application of Proposition \ref{nu-infinito}.
\qed
\vs{12}
\ni A similar result (under more restrictive conditions on $\alpha$) has been obtained by Schmid-Tretter in \cite{ScTr-02}; however, in 
\cite{ScTr-02} no information on the nodal properties of the eigenfunctions is provided.

\vs{12}
\ni
\section{The nonlinear eigenvalue problem} \label{sez-nonlineare}
\vs{12}
\ni
\subsection{A bifurcation result}
\vs{12}
\ni
In this section we are interested in proving a global bifurcation result for a nonlinear equation of the form
\begin{equation} \label{eqcompleta}
Jz'+P(x)z=\l z +S(x,z)z,\ \l\in \R,\ x>0,\ z\in \RR,
\end{equation}
where $P\in {\cal P}_{\mu}$ and $S\in C((0,+\infty)\times \RR, M^2_S)$. We denote by ${\cal S}$ the set of continuous functions $S:(0,+\infty)\times {\R}^2\lra M^{2,2}_S$ satisfying the conditions
\begin{description}
\item{$({\cal S}_1)$} there exist $\alpha \in L^{\infty}(0,+\infty)$, $\eta_{ij}\in C({\R}^2)$ such that $\eta_{ij}(0)=0$, $i, j=1, 2$, and
\begin{equation} \label{ipoS1}
\begin{array}{l}
|S_{i,j}(x,z)|\leq \alpha(x) \eta_{ij}(z),\quad \forall \ x>0,\quad z\in{\R}^2,\ i, j=1, 2;
\end{array}
\end{equation}
\ni
\item{$({\cal S}_2)$} for every compact $K\subset {\R}^2$ there exists $A_K>0$ such that
\begin{equation} \label{ipoS2}
||S(x,z)-S(x,z')||\leq A_K||z-z'||,\quad \forall \ x>0,\quad z, z'\in K.
\end{equation}
\end{description}
\vs{12}
\ni
Let $\Sigma$ denote the set of nontrivial solutions of \eqref{eqcompleta} in $D_0\times {\Lambda}$ and let $\Sigma'=\Sigma \cup \{(0,\l)\in D_0\times {\Lambda}:\ \l \ \mbox{is an eigenvalue of $A_0$}\}$, where $D_0$ and $A_0$ are as in Section \ref{sez-autovalori}. We denote by $||\cdot||_0$ the graph norm induced on $D_0$ by $A_0$, defined as
\[
||z||^2_0=||z||_{L^2(0,+\infty)}^2+||\tau z||_{L^2(0,+\infty)}^2,\quad \forall \ z\in D_0.
\]
\vs{6}
\ni
Let $M$ denote the Nemitskii operator associated to $S$, given by
\[
M(z)(x)=S(x,z(x))z(x),\quad \forall \ x>0,
\]
for every $z\in D_0$. We can show the validity of the following:
\vs{12}
\ni
\begin{proposition} \label{propM} Assume that $S\in {\cal S}$ and that
\begin{equation} \label{limitealpha}
\lim_{x\to +\infty} \alpha(x)=0,
\end{equation}
where $\alpha$ is given in \eqref{ipoS1}. Then $M:D_0\lra L^2(0,+\infty)$ is a continuous compact map and satisfies \begin{equation} \label{Mtrascurabile}
M(z)=o(||z||_0),\quad z\to 0.
\end{equation}
\end{proposition}
\vs{12}
\ni
The proof of Proposition \ref{propM} is based on the application of the following lemma:
\vs{12}
\ni
\begin{lemma} \label{convergenze} Assume that $z_0, f_0\in L^2(0,1)$ satisfy
\[
\tau z_0=f_0
\]
and let $\{z_n\}\subset L^2(0,1)$ be a sequence such that
\[
\tau z_n=f_n,
\]
for some $f_n\in L^2(0,1)$. If $z_n \rightharpoonup z_0$ and $f_n \rightharpoonup f_0$ weakly in $L^2 (0,1)$, then $Mz_n\to Mz_0$ strongly in $L^2(0,1)$.
\end{lemma}
\vs{8}
\ni
{\sc Proof.} Let us apply Theorem \ref{dicotomia-zero} to the functions $z_0$ and $z_n$, for every $n\in \N$: we have
\begin{equation} \label{decomposizione-11}
\begin{array}{l}
z_0(x)=\nu_0(x)+w_0(x),\\
\\
z_n(x)=\nu_n(x)+w_n(x),\quad \forall \ x\in (0,1),
\end{array}
\end{equation}
where
\[
\begin{array}{l}
\nu_n (x)= c_n \zeta_1 (x),\ \nu_0(x)=c_0 \zeta_1(x)\\
\\
\displaystyle{w_n(x)=\int_0^1 G(x,\xi)f_n(\xi)\,d\xi,\ w_0(x)=\int_0^1 G(x,\xi)f_0(\xi)\,d\xi,\quad \forall \ x\in (0,1).}
\end{array}
\]
Since $G\in L^\infty ((0,1)\times (0,1))$, we deduce that
\begin{equation} \label{new11}
w_n(x)\to w_0(x),\quad \forall \ x\in (0,1)
\end{equation}
by the weak convergence of $f_{n}$.
Moreover, the estimate
\begin{equation} \label{new12}
||w_n(x)- w_0(x)||\leq ||G||_{L^\infty ((0,1)^2)} ||f_n-f_0||_{L^2(0,1)},\quad \forall \ x\in (0,1),
\end{equation}
holds true; the convergence $f_n\rightharpoonup f_0$ in $L^2 (0,1)$ implies that the sequence $\{f_n\}$ is bounded in $L^2(0,1)$ and \eqref{new11}-\eqref{new12} ensure then that $w_n\to w_0$ in $L^2(0,1)$ by the dominated convergence theorem.
This condition, together with the assumption  $z_n\rightharpoonup z_0$ in $L^2 (0,1)$, implies that  $\nu_n\rightharpoonup \nu_0$ in $L^2 (0,1)$.
Hence, we obtain that $c_n\to c_0$, for $n\to +\infty$, and
\begin{equation} \label{new13}
\nu_n\to \nu_0\quad \mbox{in $L^\infty(0,1)$ and in $L^2(0,1)$.}
\end{equation}
From \eqref{new11}-\eqref{new13} we have
\begin{equation} \label{new14}
z_n(x)\to z_0(x),\quad \forall \ x\in (0,1).
\end{equation}
On the other hand, from \eqref{new12} and the boundedness of  $\{f_n\}$ in $L^2(0,1)$ we deduce also that $\{w_n\}$ is bounded in $L^\infty (0,1)$; as a consequence, using \eqref{new13}, we get that $\{z_n\}$ is bounded in $L^\infty(0,1)$. Using assumption $({\cal S}_2)$, from this fact we infer that there exists $C_1>0$ such that
\begin{equation} \label{new15}
||S(x,z_n (x))- S(x,z_0 (x))||\leq C_1 ||z_n(x)-z_0(x)||,\quad \forall \ x\in (0,1);
\end{equation}
equations \eqref{new14}-\eqref{new15} guarantee that
\begin{equation} \label{dim-88888}
S(x,z_n (x))\to S(x,z_0 (x))\quad  \forall \ x\in (0,1).
\end{equation}
Finally, from \eqref{new15} and the boundedness of  $\{z_n\}$ in $L^\infty(0,1)$ we also deduce that there exists $C_2>0$ such that
\begin{equation} \label{new16}
||S(x,z_n (x))z_{n}(x) - S(x,z_0 (x))z_{0}(x)||\leq C_2,\quad \forall \ x\in (0,1),\quad \forall \ n\geq 1;
\end{equation}
an application of the Lebesgue convergence Theorem gives
\[
\int_0^1 ||S(x,z_n (x)) z_{n}(x)- S(x,z_0 (x)) z_{0}(x)||^2\,dx\to 0,\quad n\to +\infty,
\]
i.e. $Mz_n\to Mz$ in $L^2(0,1)$.
\qed
\vs{12}
\ni
{\sc Proof of Proposition \ref{propM}.} First of all, let us observe that it is sufficient to prove the result when $x\in (0,1)$. Indeed, the fact that $P\in L^\infty(1,+\infty)$ implies that the graph norm $||\cdot||_0$, when applied to functions defined on $[1,+\infty)$, is equivalent to the $H^1(1,+\infty)$ norm; hence, when $x\in [1,+\infty)$ we can apply \cite[Prop. 4.3]{CaDa-10}.
\vs{4}
\ni
1. We first show that $Mz\in L^2(0,+\infty)$ when $z\in D_0$; from Proposition \ref{dominio} we deduce that $z\in L^\infty (0,1)$. Therefore there exists $C_z>0$ such that
\[
|S(x,z(x))|\leq C_z,\quad \forall \ x\in (0,1).
\]
As a consequence we obtain $Mz\in L^\infty(0,1)\subset  L^2(0,1)$.
\vs{6}
\ni
2. Let us fix $z_0\in D_0$ and let $z_n\in D_0$ such that $z_n\to z_0$ when $n\to +\infty$; this implies that
\begin{equation} \label{due-convergenze}
z_n\to z_0\ {\mbox{in }}\ L^2(0,1),\quad \tau z_n\to \tau z_0\ {\mbox{in }}\ L^2(0,1)
\end{equation}
We can then apply Lemma \ref{convergenze} and obtain that $Mz_n\to Mz_0$ in $L^2 (0,1)$.
\vs{6}
\ni
3. As far as the compactness of $M$ is concerned, let $\{z_n\}\subset L^2(0,1)$ be such that
\[
||z_n||_0\leq K,
\]
for some $K>0$. This implies that, up to a subsequence, we have
\begin{equation} \label{due-convergenze-2}
z_n\rightharpoonup z_0\ {\mbox{in }}\ L^2(0,1),\quad \tau z_n\rightharpoonup \tau z_0\ {\mbox{in }}\ L^2(0,1).
\end{equation}
Hence, according to Lemma \ref{convergenze}, we conclude that $Mz_n\to Mz_0$ in $L^2 (0,1)$.
\vs{6}
\ni
4. Finally, let us prove \eqref{Mtrascurabile}. We have
\begin{equation} \label{new21}
||Mz||^2_{L^2 (0,1)} =\int_0^1 ||Mz(x)||^2\,dz \leq \int_0^1 ||S(x,z(x))||^2\  ||z(x)||^2\,dx,\quad \forall \ z\in D_0.
\end{equation}
Assume now that $z\to 0$ in $D_0$; this implies that $z\to 0$ and $\tau z\to 0$ in $L^2(0,1)$; arguing as in the proof of Lemma \ref{convergenze}, we deduce that $z\to 0$ in $L^\infty (0,1)$. Therefore, assumption ${\cal S}_2$ implies that there exists $C>0$ such that
\begin{equation} \label{new22}
||S(x,z(x))|| \leq C||z(x)||\leq C ||z||_{L^\infty(0,1)},\quad \forall \ x\in (0,1).
\end{equation}
From \eqref{new21} and \eqref{new22} we deduce that
\[
||Mz||_{L^2(0,1)}\leq C ||z||_{L^\infty(0,1)} ||z||_{L^2(0,1)}\leq  C ||z||_{L^\infty(0,1)} ||z||_0,
\]
which implies that $Mz=o(||z||_0)$ as $||z||_0\to 0$.
\qed
\vs{12}
\ni
Now, let us observe that, in view of the results on $A_0$ given in Section \ref{sez-autovalori} and of Proposition \ref{propM}, it is possible to write \eqref{eqcompleta} as an abstract equation of the form
\begin{equation} \label{astratta}
A_0 u+M(u)=\l u, \ (u,\l)\in D_0\times {\R},
\end{equation}
where $A_0:D_0\subset L^2 (0,+\infty)\to L^2 (0,+\infty)$ is an unbounded self-adjoint operator such that
\[
\sigma_{ess} (A_0)=(-\infty,\mu^-]\cup [\mu^+,+\infty)
\]
and $M:D_0\times {\R}\lra L^2(0,+\infty)$ is a continuous and compact map such that
\begin{equation} \label{Mtrascurabile-astratto}
M(u)=o(||u||),\quad u\to 0.
\end{equation}
From an application of a global bifurcation result (see \cite[Th. 1.2]{St-75}, \cite[Th. 4.1]{CaDa-10}) to \eqref{astratta} we then obtain the following main result:
\vs{12}
\ni
\begin{Theorem} \label{main-alternative} Assume that $P\in {\cal P}_{\mu}$, $S\in {\cal S}$ and that \eqref{limitealpha} holds true. Then, for every eigenvalue $\gamma\in \Lambda$ of $A_0$ there exists a continuum $C_\gamma$ of nontrivial solutions of \eqref{eqcompleta} in $D_0\times {\R}$ bifurcating from $(0,\gamma)$ and such that one of the following conditions holds true:
\begin{description}
\item{(1)} $C_\gamma$ is unbounded in $D_0\times {\Lambda}$;
\item{(2)} $\sup \{\l:\ (u,\l)\in C_\gamma\}\geq \mu^+$ or $\inf \{\l:\ (u,\l)\in C_\gamma\}\leq \mu^-$;
\item{(3)} $C_\gamma$ contains $(0,\gamma')\in \Sigma'$, with $\gamma'\neq \gamma$.
\end{description}
\end{Theorem}
\vs{12}
\ni
Now, let us observe that a more precise description of the bifurcating branch, eventually leading to exclude condition (3),  can be obtained when there exists a continuous functional $i:\Sigma'\to \Z$ (cf. \cite[Th. 4.2]{CaDa-10}). In order to define such a functional, we first define the rotation number of solutions to \eqref{eqcompleta} by means of a linearization procedure; to this aim for every solution $(w,\mu)$ of \eqref{eqcompleta} we consider the linear equation
\begin{equation} \label{linearizzato}
                         Jz'+P(x)z=\mu z+S(x,w(x))z,
\end{equation}
which obviously reduces to
\begin{equation} \label{linearizzato2}
                         Jz'+P(x)z=\mu z
\end{equation}
when $w=0$. It is clear that $w$ is a solution of \eqref{linearizzato}; let us denote by $P_w$ the matrix defined by
\[
P_w(x)=P(x)-S(x,w(x)),\quad \forall \ x>0.
\]
We can prove the following result:
\vs{12}
\ni
\begin{lemma} \label{prop-linearizzato} For every $(w,\mu)\in \Sigma$ we have $P_w\in {\cal P}_{\mu}$.
\end{lemma}
\vs{8}
\ni
{\sc Proof.} Let us first observe that $w\in D_0$ implies that $w\in H^1 (1,+\infty)$ and $w\in L^\infty (0,+\infty)$ (cf. Proposition \ref{dominio}). In particular we have
\begin{equation} \label{dim-011}
\lim_{x\to +\infty} w(x)=0
\end{equation}
and  there exists a compact set $K_w\subset \RR$ such that
\[
w(x)\in K_w,\quad \forall \ x>0.
\]
Using \eqref{limitiPzero}, assumption $({\cal S}_1)$ and the fact that $w\in L^\infty (0,1)$, we obtain that
\[
\lim_{x\to 0^+} x^{\beta} P_w(x)=\lim_{x\to 0^+} ( x^{\beta} P(x)-x^{\beta} S(x,w(x)))=P^*;
\]
therefore $P_w$ satisfies \eqref{limitiPzero}. Moreover, we have
\begin{equation} \label{dim33}
R_{0,w}(x)= x^{\beta} P_w(x)-P^*=R_0(x)- x^{\beta} S(x,w(x)),\quad \forall \ x>0.
\end{equation}
Using again $({\cal S}_1)$,  we plainly deduce that there exists $\eta \in C(\RR,\R^+)$ such that
\begin{equation} \label{dim44}
\int_0^{1} \dfrac{1}{x^{\beta} } ||x^{\beta} S(x,w(x))||^{q_0}\,dx= \int_0^{1}  x^{\beta(q_0-1)}  \alpha (x)\eta (w(x))\,dx<+\infty,
\end{equation}
since $q_0\geq 1$ and $\alpha, w\in L^\infty (0,1)$. From \eqref{iporestozero}, \eqref{dim33} and \eqref{dim44} we can conclude that $R_{0,w}$ satisfies \eqref{iporestozero}.
\vs{6}
\ni
Now, we pass to the proof of the validity of $({\cal P}_1)$. Using \eqref{limitiPinfinito}, assumption $({\cal S}_1)$ and \eqref{dim-011}, we infer that
\[
\lim_{x\to +\infty} P_w(x)=\lim_{x\to +\infty} (P(x)-S(x,w(x)))=P_\infty;
\]
hence $P_w$ satisfies \eqref{limitiPinfinito}.
\vs{4}
\ni
Moreover, we have
\begin{equation} \label{dim11}
R_{\infty,w}(x)=P_w(x)-P_\infty=R_\infty(x)-S(x,w(x)),\quad \forall \ x>0.
\end{equation}
From assumption $({\cal S}_2)$, with $K=K_w$ and $z'=0$, we obtain
\begin{equation} \label{dim22}
\int_1^{+\infty} ||S(x,w(x))||^2\,dx\leq A_{K_w}^2 \int_1^{+\infty} ||w(x)||^2\,dx<+\infty.
\end{equation}
When $q_\infty\geq 2$ this allows to conclude that $R_{\infty,w}$ satisfies \eqref{iporestoinfinito}, since
\begin{equation} \label{dim33a}
\begin{array}{l}
\displaystyle{\int_1^{+\infty} ||S(x,w(x))||^{q_\infty}\,dx= \int_1^{+\infty} ||S(x,w(x))||^{q_\infty-2}\ ||S(x,w(x))||^{2}\,dx\leq} \\
\\
\displaystyle{\leq C_w \int_1^{+\infty} ||S(x,w(x))||^2\,dx<+\infty.}
\end{array}
\end{equation}
Finally, also when $q_\infty<2$ it is possible to show that $R_{\infty,w}$ satisfies \eqref{iporestoinfinito} with the same ${q_\infty}$ of $R_\infty$; indeed, at this point we can say that $w\in H^1 (1,+\infty)$ is a nontrivial solution of the linear equation
\[
Jz'+P_w (x) z=\mu z,
\]
where $P_w\in {\cal P}_{\mu}$ and $\mu \in \Lambda$. Therefore Proposition \ref{stimeinf2} applies and we deduce that $w$ satisfies the first condition in \eqref{stimez}. As a consequence, $w\in L^{q_\infty} (1,+\infty)$ and we are able to repeat \eqref{dim22} with the exponent $q_\infty$ instead of $2$.
\qed
\vs{12}
\ni
As a consequence of Lemma \ref{prop-linearizzato}, the results of Section \ref{sistemalineare} apply to \eqref{linearizzato}; in particular, when $w\neq 0$ we can consider the number $\rot \, (w)$ defined in \eqref{defrotazioni}.
\vs{12}
\ni
\begin{Definition} \label{def-rot-2} Assume that $P\in {\cal P}_{\mu}$ and $S\in {\cal S}$ and let $(w,\mu)$ be a solution of \eqref{eqcompleta}.
\vs{4}
\ni
If $(w,\mu)\neq (0,\mu)$, then the rotation number $j(w,\mu)$ of $(w,\mu)$ is defined by
\begin{equation} \label{def-j-1}
j(w,\mu)=\rot \, (w).
\end{equation}
If $(w,\mu)=(0,\mu)$ and the linear problem \eqref{linearizzato2}
has a nontrivial solution $z_\mu$ belonging to $H^1(0,+\infty)$, then the rotation number $j(w,\mu)$ of $(w,\mu)$  is defined
by
\begin{equation} \label{def-j-2}
j(w,\mu)=\rot \, (z_\mu).
\end{equation}
\end{Definition}
\vs{12}
\ni
By means of Definition \ref{def-rot-2} we have defined $j:\Sigma'\to \R$; this functional will be used in order to construct a continuous discrete functional whose values are preserved in the bifurcating branches $C_\gamma$ of solutions of \eqref{eqcompleta}. It is important now to observe that every branch $C_\gamma$ satisfies
\[
C_\gamma \subset H^1_0(0,+\infty)\times \R;
\]
indeed, this is a consequence that $(z,\l)\in C_\gamma$ is a solution of the linear equation
\[
Jz'+P_z(x)z=\l z
\]
such that $z\in D_0$. According to Remark \ref{regolarita-sol} this implies that $z\in H^1_0 (0,+\infty)$.
\vs{6}
\ni
Hence, it is sufficient to study the continuity properties of $j$ with respect to the $H^1_0 (0,+\infty)$-norm, denoted by $|||\cdot |||$.
\vs{12}
\ni
\begin{proposition} \label{continuita} The function $j:\Sigma'\to \R$ is continuous.
\end{proposition}
\vs{8}
\ni
{\sc Proof.} We prove the continuity of $j$ at every point $(w,\mu)\in \Sigma \cap H^1_0 (0,+\infty)$. In a very similar way it is possible to show that $j$ is also continuous at every point $(0,\l)$, with $\l$ eigenvalue of $A_0$.
\vs{6}
\ni
Let us fix  $(w,\mu)\in \Sigma \cap H^1_0 (0,+\infty) $ and let $\epsilon >0$ small enough; consider then the numbers $\delta$, $\delta_0$, $x_\infty$ and $x_0$ given in Proposition \ref{confini-theta-inf} and Proposition \ref{confini-theta-zero} (with $\tilde \l=\mu$ and $\tilde P=P_w$) and let $\delta_1=\min (\delta,\delta_0)$.
\vs{4}
\ni
Using assumption $({\cal S}_1)$ and the continuous embedding $H^1_0 (0,+\infty)\subset L^\infty (0,+\infty)$, it is possible to show that there exist $\delta_2 >0$ such that
\[
||P_z-P_w||_{L^\infty (0,+\infty)}<\delta_1
\]
if $|||z-w|||<\delta_2$.
\vs{4}
\ni
Hence, from Proposition \ref{confini-theta-inf} and Proposition \ref{confini-theta-zero} we deduce that for every $(z,\l)\in \Sigma \cap H^1_0 (0,+\infty) $ with $|\l-\mu|<\delta_2$ and $|||z-w|||<\delta_2$ we have
\begin{equation} \label{confini-theta-nonlineare}
\begin{array}{l}
|\theta_z (x,\l)-\theta_z (+\infty,\l)|<\epsilon,\quad \forall \ x\geq x_{\infty}\\
\\
|\theta_z (x,\l)-\theta_z (0)|<\epsilon,\quad \forall \ x\in (0,x_{0}].
\end{array}
\end{equation}
Now, let us observe that we have
\[
\begin{array}{l}
j(z,\l)-j(w,\mu)=\dfrac{\theta_z (+\infty,\l)-\theta_z (0)}{\pi}-\dfrac{\theta_w (+\infty,\mu)-\theta_w (0)}{\pi}=\dfrac{\theta_z (+\infty,\l)-\theta_w (\infty,\mu)}{\pi},
\end{array}
\]
since $\theta_z(0)=\theta_w(0)$. Therefore, the result follows from the same argument used in the proof of Proposition \ref{continuita-nu*}.
\qed
\vs{12}
\ni
Before defining the functional $i$ a remark is in order; we recall that if $(w,\l)\in \Sigma \cap H^1_0 (0,+\infty)$ then
\[
\theta_w (0)=\arctan \dfrac{w^*_{1,2}}{w^*_{1,1}}\in (0,\pi)  \quad {\mbox{mod }} \, \pi
\]
and
\[
\theta_w (+\infty,\l)=\pi-\arctan \sqrt{\dfrac{\l-\mu^-}{\mu^+-\l}}\in \left(\dfrac{\pi}{2},\pi\right)\quad {\mbox{mod }} \, \pi.
\]
As a consequence, when
\begin{equation} \label{primoquadrante2}
\arctan \dfrac{w^*_{1,2}}{w^*_{1,1}}\in \left(0,\dfrac{\pi}{2}\right)
\end{equation}
we have
\[
j(w,\l)\not\in \Z,\quad \forall \ (w,\l)\in \Sigma \cap H^1_0 (0,+\infty).
\]
On the other hand, if
\begin{equation} \label{secondoquadrante2}
\arctan \dfrac{w^*_{1,2}}{w^*_{1,1}}\in \left(\dfrac{\pi}{2},\pi \right)
\end{equation}
we have
\[
j(w,\l)+\dfrac{1}{2}\not\in \Z,\quad \forall \ (w,\l)\in \Sigma \cap H^1_0 (0,+\infty).
\]
This suggests to define $i: \Sigma'\to \Z$ as
\begin{equation} \label{def-i-1}
i(w,\l)=\left[j(w,\l)\right],\quad \forall \ (w,\l)\in \Sigma',
\end{equation}
if \eqref{primoquadrante2} holds true, and
\begin{equation} \label{def-i-2}
i(w,\l)=\left[j(w,\l)+\dfrac{1}{2}\right],\quad \forall \ (w,\l)\in \Sigma',
\end{equation}
if \eqref{secondoquadrante2} holds true (recall also Proposition \ref{numero-rotazioni-autofunzioni}). Let us observe that Proposition \ref{numero-rotazioni-autofunzioni} also implies that
\begin{equation} \label{indici-diversi-due}
i(z_\gamma,0)\neq i(z_{\gamma'},0),
\end{equation}
for every $\gamma \neq \gamma'\in \Lambda$ eigenvalues of $A_0$ (with associated eigenfunctions $z_\gamma$ and $z_{\gamma'}$, respectively).
\vs{4}
\ni
From Proposition \ref{continuita} and the definition of $i$ we obtain the following result:
\vs{12}
\ni
\begin{proposition} \label{continuita-i} The function $i:\Sigma'\to \R$ is continuous.
\end{proposition}
\vs{12}
\ni
As a consequence, using Proposition \ref{continuita-i} and \eqref{indici-diversi-due}, from Theorem \ref{main-alternative} we deduce the final result:
\vs{12}
\ni
\begin{Theorem} \label{main-alternative-2} Assume that $P\in {\cal P}_{\mu}$, $S\in {\cal S}$ and that \eqref{limitealpha} hold true. Then, for every eigenvalue $\gamma\in \Lambda$ of $A_0$ there exists a continuum $C_\gamma$ of nontrivial solutions of \eqref{eqcompleta} in $D_0\times {\R}$ bifurcating from $(0,\gamma)$ and such that one of the conditions (1)-(2) of Theorem \ref{main-alternative} holds true and
\begin{equation} \label{i-costante}
i(w,\l)=i(z_\gamma,0),\quad \forall \ (w,\l)\in C_\gamma,
\end{equation}
where $z_\gamma$ is the eigenfunction of $A_0$ associated to $\gamma$.
\end{Theorem}
\vs{12}
\ni
\subsection{Application to the Dirac equation} \label{Dirac-nonlineare}
\vs{12}
\ni
Let us consider the partial differential equation
\begin{equation} \label{Dirac-1}
i\sum_{j=1}^3 \alpha_j \dfrac{\partial \psi}{\partial x_j}-\beta \psi - V(||x||)\psi + i a \sum_{j=1}^3 \alpha_j \dfrac{\partial V(||x||)}{\partial x_j} \psi  =\l \psi+\gamma (||x||)F(\langle\beta\psi,\psi\rangle)\beta\psi,\quad x\in {\R}^3, \ a\in \R,
\end{equation}
where $\psi:{\R}^3\to {\C}^4$, $V\in C((0,+\infty),\R)$ satisfies \eqref{ipotesi-potenziale-infinito}-\eqref{ipotesi-potenziale-zero}-\eqref{ipotesi-potenziale-zero-resto1}-\eqref{ipotesi-potenziale-zero-resto2}, $\gamma \in C((0,+\infty),\R)$ fulfills
\begin{equation} \label{ipotesi-gamma}
\lim_{r\to 0^+} r^2\gamma (r)\in \R,\quad r^2\gamma (r)=o(1),\ r\to +\infty,
\end{equation}
$F\in C(\R,\R)$, $\langle \cdot,\cdot \rangle$ denotes the scalar product in ${\C}^4$ and $\alpha_j$ ($j=1, 2, 3$) and $\beta$ are the $4\times 4$ matrices given by
\[
\alpha_j=\left(\begin{array}{cc}
                0&\sigma_j\\
                &\\
                \sigma_j&0
                \end{array}
         \right),\quad
\beta=\left(\begin{array}{cc}
                \sigma_0&0\\
                &\\
                0&-\sigma_0
                \end{array}
         \right),
\]
where
\[
\sigma_0=\left(\begin{array}{cc}
                1&0\\
                &\\
                0&1
                \end{array}
         \right),\quad
\sigma_1=\left(\begin{array}{cc}
                0&1\\
                &\\
                1&0
                \end{array}
         \right),\quad
\sigma_2=\left(\begin{array}{cc}
                0&-i\\
                &\\
                i&0
                \end{array}
         \right),\quad
\sigma_3=\left(\begin{array}{cc}
                1&0\\
                &\\
                0&-1
                \end{array}
         \right).
\]
We remark that nonlinearities like the one in \eqref{Dirac-1} give rise to the so-called generalized Soler models (see \cite{Es-02}).
In fact, Soler \cite{So-70} formulated a model of extended fermions by introducing a self interaction term which corresponds to the choice $ F(s) = s $ in \eqref{Dirac-1}
(see \cite{Ra-book-83} for a survey on interaction terms which are interesting from a physical point of view).

We denote by $H_0$ the (free) Dirac operator defined by
\begin{equation} \label{def-H0}
H_0\psi=i\sum_{j=1}^3 \alpha_j \dfrac{\partial \psi}{\partial x_j}-\beta \psi,\quad \forall \ \psi \in H^1_0({\R}^3)\subset L^2({\R}^3).
\end{equation}
In \cite{Th-book-92} a decomposition of $H_0 - V+ia\ \alpha\cdot \nabla V $ has been performed, using polar coordinates in ${\R}^3$ and the unitary isomorphism
\begin{equation} \label{isomorfismo}
\begin{array}{ll}
\varphi:&L^2({\R}^3)\to L^2((0,+\infty),dr;L^2(S^2))\\
&\\
&\psi\mapsto \tilde \psi,
\end{array}
\end{equation}
where $\tilde \psi$ is defined by
\begin{equation} \label{def-iso}
\tilde \psi (r,\theta,\phi)=r\psi (x(r,\theta,\phi)),\quad \forall \ r>0,\ (\theta,\phi)\in S^2.
\end{equation}
In order to describe such a decomposition, for every $l=0, 1, 2, \ldots$ and $m=-l, -l+1, \ldots, l$ let us denote by $Y^m_l$ the usual spherical harmonic; moreover, for every $j=1/2, 3/2, 5/2, \ldots$, let $m_j=-j, -j+1,\ldots, j$ and $k_j=-(j+1/2), j+1/2$ and define
\begin{equation} \label{def-Psi}
\Psi^{m_j}_{j-1/2} =\dfrac{1}{\sqrt{2j}} \left(
                                                                    \begin{array}{l}
                                                                    \sqrt{j+m_j}\  Y^{m_j-1/2}_{j-1/2}\\
                                                                    \\
                                                                    \sqrt{j-m_j}\ Y^{m_j+1/2}_{j-1/2}
						    \end{array}
						    \right),\quad
\Psi^{m_j}_{j+1/2} =\dfrac{1}{\sqrt{2j+2}} \left(
                                                                    \begin{array}{l}
                                                                    \sqrt{j+1-m_j}\  Y^{m_j-1/2}_{j+1/2}\\
                                                                    \\
                                                                   - \sqrt{j+1+m_j}\ Y^{m_j+1/2}_{j+1/2}
						    \end{array}
						    \right)
\end{equation}
and
\begin{equation} \label{def-Phi}
\Phi^{+}_{m_j,\mp (j+1/2)} =\left(
                                                                    \begin{array}{l}
                                                                   i\ \Psi^{m_j}_{j\mp 1/2}\\
                                                                    \\
                                                                    0
						    \end{array}
						    \right),\quad
\Phi^{-}_{m_j,\mp (j+1/2)}  =\left(
                                                                    \begin{array}{l}
                                                                   0\\
                                                                    \\
                                                                   \Psi^{m_j}_{j\pm 1/2}
						    \end{array}
						    \right).
\end{equation}
We also set
\begin{equation} \label{def-spazi}
{{\cal H}}_{m_j,k_j}={\mbox{span }} (\Phi^+_{m_j,k_j},\Phi^-_{m_j,k_j}),\quad \forall \ j=1/2, 3/2, \ldots.
\end{equation}
Then, we have the following result:
\vs{12}
\ni
\begin{Theorem} \label{decomposizione} (\cite[Th. 4.14]{Th-book-92}) For every $j=1/2, 3/2, \ldots$ the subspace 
$C_0^\infty (0,+\infty)\otimes {\cal H}_{m_j,k_j}\subset L^2((0,+\infty),dr;L^2(S^2){^4})$ is invariant under the action of $H_0 -V+ia\ \alpha\cdot \nabla V$.
Moreover, with respect to the basis $\{\Phi^+_{m_j,k_j}, \Phi^-_{m_j,k_j}\}$ of ${\cal H}_{m_j,k_j}$ the restriction of $H_0 -V+ia\ \alpha\cdot \nabla V $ to ${\cal H}_{m_j,k_j}$ can be represented by the operator ${ h}_{m_j, k_j}$ given by
\begin{equation} \label{def-Dirac-radiale}
{h}_{m_j,k_j}=\left(\begin{array}{cc}
                                        \displaystyle{-1 \ -V \  }&{\displaystyle{\dfrac{d}{dr}-\dfrac{k_j}{r} +aV'}}\\
 			        &\\
			         {\displaystyle{-\dfrac{d}{dr}-\dfrac{k_j}{r} +aV'}}& 1 -V
				\end{array}
			\right).
\end{equation}
Moreover, the Dirac operator $H_0 -V+ia\ \alpha\cdot \nabla V$ on $C_0^{\infty} ({\R}^3)^4$ is unitarily equivalent to the direct sum of the partial waves operators ${ h}_{m_j,k_j}$, i.e.
\[
\begin{array}{l}
H_0 -V+ia\ \alpha\cdot \nabla V\approx \bigoplus_{j=1/2, 3/2,\ldots}^{+\infty} \ \bigoplus_{m_j=-j}^j\ \bigoplus_{k_j=\pm (j+1/2)} {{h}}_{m_j,k_j}.
\end{array}
\]
\end{Theorem}
\vs{12}
\ni

\begin{remark} \label{significato}
The partial wave subspaces can be considered as a suitable generalization of radial functions adapted to the structure of the nonlinear problem. More precisely, the vectors $\Phi^{\pm}_{m_j,k_j}$ which are a basis for the partial wave subspace ${{\cal H}}_{m_j,k_j}$ are the eigenfunctions of the spin orbit operator (cf. (\cite{Th-book-92}). We also observe that these subspaces are implicitly used in \cite{BCDM-88},\cite{BaCaVa-90}, where (having in mind the Soler model) the system of ODEs is obtained from the PDE by making the ansatz that solutions should be a linear combination of functions of the form ${\Phi}^{+}_{1/2,1}=(\dfrac{i}{2\sqrt{\pi}} \cos \theta,\dfrac{i}{2\sqrt{\pi}} e^{i \phi}\sin \theta,0,0), {\Phi}^{-}_{1/2,1}=(0,0,\dfrac{1}{2\sqrt{\pi}},0)$. On the same lines but in the context of the Schr\"odinger equation, we refer to \cite[Example 1.5]{We-book-87}.
\end{remark}

\vs{12}
\ni
Let us observe that the operators $\tau_{k_j}={h}_{m_j, k_j}$, $j=1/2, 3/2, \ldots$, are of the form \eqref{tau} with $P=P_{V,k_j,a} $ as in \eqref{Dirac-magn}. Therefore, we can apply the theory developed in Sections \ref{sistemalineare} and \ref{sez-autovalori}; in particular we can consider the selfadjoint realization $A_0$ of $\tau_{k_j}$, $j=1/2, 3/2, \ldots$, defined in \eqref{def-A0}. We denote by $A_{k_j}$ this operator, by $D_{k_j}$ its domain and 
we define
\[
 E_{k_j} = \{ u^{+} \Phi^{+}_{m_j,k_j} + u^{-} \Phi^{-}_{m_j,k_{j}} : u = ( u^{+}, u^{-} ) \in D_{ k_{j} } \}.
\]

\ni
From Theorem \ref{decomposizione} and the definition of $D_{k_j}$ we immediately deduce that the image of $E_{k_j}$ via the operator $H_0 -V+ia\ \alpha\cdot \nabla V$ is contained in $L^{2}(( 0, \infty ))\otimes {\cal H}_{ k_{j}, m_{j} },$ for every $j=1/2, 3/2, \ldots$.
\vs{6}
\ni
Now, let us observe that Theorem \ref{decomposizione} states that the subspaces
\[
C_0^\infty (0,+\infty)\otimes {\cal H}_{m_j,k_j}\subset L^2((0,+\infty),dr;L^2(S^2)^4), \quad j=1/2, 3/2, \ldots,
\]
are preserved by the linear operator in \eqref{Dirac-1}. It is important to note that in the particular case of $j=1/2$ the subspaces 
$C_0^\infty (0,+\infty)\otimes {\cal H}_{m_{1/2},k_{1/2}}\subset L^2((0,+\infty),dr;L^2(S^2)^4)$
are invariant also for the nonlinear term $F(\langle\beta\psi,\psi\rangle)\beta\psi$ in \eqref{Dirac-1} (cf. \cite[Lemma 5.5]{Ca-11}), when $F$ is regular.
\vs{4}
\ni
Indeed, let $u\in L^2({\R}^3)^4$ such that $\varphi (u)\in C_0^\infty (0,+\infty)\otimes {\cal H}_{m_{1/2},k_{1/2}}$, where $\varphi$ is defined in \eqref{isomorfismo}-\eqref{def-iso}. A simple computation, based on the expressions of the functions $\Phi^\pm _{m_{1/2},k_{1/2}}$,
shows that, if we have
\[
(\varphi(u))(r,\theta,\phi)=u^+(r)\Phi^+_{m_{1/2},k_{1/2}}(\theta,\phi)+u^-(r)\Phi^-_{m_{1/2},k_{1/2}}(\theta,\phi), 
\]
then
\begin{equation} \label{conto-nonlinearita-1}
\langle \beta u(x),u(x) \rangle =\dfrac{1}{4\pi^2 r^2} [(u^+(r))^{2} - (u^-(r))^2]
\end{equation}
and
\begin{equation} \label{conto-nonlinearita-2}
\begin{aligned}
F & \left( \left\langle \beta \dfrac{(\varphi(u))(r,\theta,\phi)}{r},\dfrac{(\varphi(u))(r,\theta,\phi)}{r} \right\rangle \right) \beta[(\varphi(u))(r,\theta,\phi)] \\
& = F\left( \dfrac{1}{4\pi^2 r^2} (u^+(r))^{2} - (u^-(r))^2\right) \left( u^+(r) \Phi^+_{m_{1/2},k_{1/2}}(\theta,\phi) - u^-(r)\Phi^-_{m_{1/2},k_{1/2}}(\theta,\phi) \right),
\end{aligned}
\end{equation}
proving that
\begin{equation} \label{conto-nonlinearita-3}
\tilde\psi\in C_0^\infty (0,+\infty)\otimes {\cal H}_{m_{1/2},k_{1/2}}\quad \Rightarrow \quad F\left(\left\langle \beta\dfrac{\tilde\psi}{r},\dfrac{\tilde\psi}{r} \right\rangle \right)\beta\tilde\psi\in C_0^\infty (0,+\infty)\otimes {\cal H}_{m_{1/2},k_{1/2}}.
\end{equation}
Then, with an argument similar to the one developed in the proof of Proposition \ref{propM}, we deduce that
\begin{equation} \label{conto-nonlinearita-4}
\tilde\psi\in E_{k_{1/2}} \quad \Rightarrow \quad \gamma (r) F\left(\left\langle \beta \dfrac{\tilde\psi}{r},\dfrac{\tilde\psi}{r} \right\rangle \right)\beta\tilde\psi\in L^2(0,+\infty)\otimes {\cal H}_{m_{1/2},k_{1/2}}.
\end{equation}
\noindent
This fact is important to obtain a relation between solutions of \eqref{Dirac-1} and solutions of a nonlinear ordinary differential equation of the form \eqref{eqcompleta}. 
Indeed, for every function $u\in L^2({\R}^3)^4$ with $\varphi (u)\in E_{k_{1/2}}$, let $z=(u^+,u^-)\in  D_{k_{1/2}}$ such that
\begin{equation} \label{relazione-soluzioni}
\varphi (u)=u^+\Phi^+_{m_{1/2},k_{1/2}}+u^-\Phi^-_{m_{1/2},k_{1/2}}.
\end{equation}
Then, \eqref{conto-nonlinearita-4} implies that $u\in L^2({\R}^3)^4$ with $\varphi (u)\in  E_{k_{1/2}}$ is a nontrivial solution of \eqref{Dirac-1} if and only if $z=(u^+,u^-)\in  D_{k_{1/2}}$ is a nontrivial solution of
\begin{equation} \label{eqcompleta-Dirac}
\tau_{k_{1/2}}z=\l z + \gamma(r) F\left( \dfrac{ ( u^{+} )^{2} - ( u^{-} )^{2} }{ 4 \pi r^2 } \right)
\left[
\begin{array}{rr}
1 &  0 \\
0 & -1
\end{array}
\right] z, \quad r>0.
\end{equation}
Let us denote $E=\varphi^{-1} (E_{k_{1/2}})$; 
in view of the above arguments and choosing
\[
S(r,z) = \gamma(r) F\left( \dfrac{ ( u^{+} )^{2} - ( u^{-} )^{2} }{ 4 \pi r^2 } \right)
\left[
\begin{array}{rr}
1 &  0 \\
0 & -1
\end{array}
\right],
\]
from Theorem \ref{main-alternative-2} we plainly obtain the following result:
\vs{12}
\ni
\begin{Theorem} \label{bif-Dirac}
Let us suppose that $V\in C(0,+\infty)$ and $\gamma \in C(0,+\infty)$ satisfy  \eqref{ipotesi-potenziale-infinito}-\eqref{ipotesi-potenziale-zero}-\eqref{ipotesi-potenziale-zero-resto1}
and \eqref{ipotesi-gamma} and let $ F : \RR \to \RR $ be a locally Lipschitz continuous function such that $ | F(s) | \le C |s| $ for all $ s \in \RR $ and some constant $ C > 0 $.
Then, for every eigenvalue $\mu\in (-1,1)$ of $A_{k_{1/2}}$ there exists a continuum $C_\mu$ of nontrivial solutions of \eqref{Dirac-1} in $E\times {\R}$ such that one of the conditions
\begin{description}
\item{(1)} $C_\mu$ is unbounded in $E\times (-1,1)$
\item{(2)} $\sup \{\l:\ (u,\l)\in C_\mu\}\geq 1$ or $\inf \{\l:\ (u,\l)\in C_\mu\}\leq -1$
\end{description}
holds true and
\begin{equation} \label{i-costante-2}
{\tilde i}(w,\l)=i(z_\mu,0),\quad \forall \ (w,\l)\in C_\mu,
\end{equation}
where
\[
{\tilde i} (w,\l)=i((w^+,w^-),\l)
\]
and $z_\mu$ is the eigenfunction of $A_{k_{1/2}}$ associated to $\mu$.
\end{Theorem}
\vs{12}
\ni

\end{document}